\documentclass[12pt]{article}

\newcommand{\eqd}{\stackrel{\rm{def}}{=}}
\newcommand{\by}{\overline{\cal Y}}
\newcommand{\br}{\overline{R}}
\usepackage{amsmath}
\usepackage{graphicx}
\begin{document}

\title{Bartlett correction of an independence test in a multivariate Poisson model}
\author{Rolf Larsson\footnote{Department of Mathematics, Uppsala University, P.O.Box 480, SE-751 06 Uppsala, Sweden. e-mail \texttt{rolf.larsson@math.uu.se}}}
\maketitle

\begin{abstract}
We consider a system of dependent Poisson variables, where each variable is the sum of an independent variate and a common variate. It is the common variate that creates the dependence. 

Within this system, a test of independence may be constructed where the null hypothesis is that the common variate is identically zero. In the present paper, we consider the maximum log likelihood ratio test. For this test, it is well-known that the asymptotic distribution of the test statistic is an equal mixture of zero and a chi square distribution with one degree of freedom. We examine a Bartlett correction of this test, in the hope that we will get better approximation of the nominal size for moderately large sample sizes. 

This correction is explicitly derived, and its usefulness is explored in a simulation study. For practical purposes, the correction is found to be useful in dimension two, but not in higher dimensions.
\end{abstract}

\section{Introduction}

Bartlett correction is a simple technique to obtain small sample corrections of asymptotic distributions, in particular useful for maximum log likelihood ratio tests. It was introduced in Bartlett (1937). Let the sample size be $n$, and write $Q_n=-2\log\Lambda_n$, where $\Lambda_n$ is the ratio between the maximum likelihoods under the null and the alternative hypotheses, respectively.
Denoting the distributional limit of $Q_n$ by $Q_\infty$, the simple idea is to write 
$$Q_n\approx \frac{E Q_n}{E Q_\infty} Q_\infty,$$
where 
$$EQ_n= E Q_\infty+n^{-1}R+o(n^{-1}),$$
where $E$ is the expectation operator and $o(n^{-1})$ is an error term that tends to zero faster than $n^{-1}$ as $n\to\infty$. Hence, 
\begin{equation}
Q_n\approx(1+n^{-1}R')Q_\infty,\label{Bartlett0}
\end{equation}
for a constant $R'=R/E Q_\infty$. Obviously, the approximation in (\ref{Bartlett0}) asymptotically corrects the first moment of the statistic, in the sense that the error is of smaller order than $n^{-1}$. In fact, for maximum log likelihood ratio tests, Lawley (1956) shows that under regularity conditions, all higher moments are corrected in the same way. Hayakawa (1977) went on to prove that this practice also corrects the fractiles of the distribution.

In standard situations, it is well known that, under the null hypothesis, the asymptotic distribution of the maximum log likelihood ratio (LR) test ($Q_\infty$ above) is $\chi^2$ with one degree of freedom, $\chi^2(1)$. However, the present paper deals with a case where the parameter of the null hypothesis is on the boundary of the parameter space. Here, the standard result does not apply. Indeed, Self and Liang (1987) have shown that the corresponding distribution is an equal mixture of a point mass at one and $\chi^2(1)$. If the parameter is $\lambda$ (say), and we test null hypothesis $H_0$: $\lambda=0$ vs $H_1$: $\lambda>0$, then the asymptotic probability that the maximum likelihood estimator $\hat\lambda=0$ is $1/2$, and the asymptotic distribution of $(Q_n|\hat\lambda>0)$ is $\chi^2(1)$.

In the present paper, we will deal with a dependent Poisson model, see e.g. Karlis (2003). This model is as follows: Assume that $U$, $X_1,X_2,...,X_m$ are independent Poisson with parameters $\lambda$, $\mu_1,\mu_2,...,\mu_m$, respectively. Then, construct the  dependent Poisson variates $(Y_1,Y_2,...,Y_m)$ through $Y_j=U+X_j$, $j=1,2,...,m$. It is clear that the $Y_j$ are simultaneously independent if and only if $\lambda=0$. Hence, it is of great interest to test $H_0$: $\lambda=0$ vs $H_1$: $\lambda>0$. Doing this with the LR test is the focus of the present paper. 

One application of such a test can be found in Larsson (2020), where it can be seen as a building block in the search algorithm for the best discrete factor analysis model. (There, the best model is defined as the one with the smallest AIC, but comparing AIC values is equivalent to using the LR test on some significance level, if the null model is a special case of the alternative model.)

The rest of the paper unfolds as follows. Section 2 gives the asymptotic approximation of the expectation of the LR test under the null. In section 3, Bartlett correction is discussed. Section 4 gives small sample simulation results, and section 5 concludes.

\section{Asymptotic approximation of the expectation}

Suppose we have a random sample $(y_{i1},y_{i2},...,y_{im})$, $i=1,2,...,n$, from the distribution of $(Y_1,Y_2,...,Y_m)$ of the dependent Poisson model by Karlis (2003).
We want to test $H_0$: $\lambda=0$ (independence) vs $H_1$: $\lambda>0$ (dependence) by using the likelihood ratio test.

The first object is to find the maximum likelihood estimates (MLEs) of the parameters. 

Let $y_i=\min(y_{i1},y_{i2},...,y_{im})$ for all $i$. The likelihood is

\begin{align}
&L(\lambda,\mu_1,\mu_2,...,\mu_m)\notag\\&=\prod_{i=1}^n \sum_{u=0}^{y_i}e^{-\lambda}\frac{\lambda^u}{u!}\prod_{j=1}^m e^{-\mu_j}\frac{\mu_j^{y_{ij}-u}}{(y_{ij}-u)!}\notag\\
&=e^{-n(\lambda+\sum_{j=1}^m\mu_j)}\left(\prod_{j=1}^m\mu_j^{n\bar y_j}\right)\prod_{i=1}^n\sum_{u=0}^{y_i}\frac{1}{u!\prod_{j=1}^m(y_{ij}-u)!}\left(\frac{\lambda}{\prod_{j=1}^m\mu_j}\right)^u,\label{L}
\end{align}
where $\bar y_j=n^{-1}\sum_{i=1}^n y_{ij}$ for $j=1,2,...,m$.

Let $\hat\lambda$ and $\hat\mu_j$ be the MLEs of $\lambda$ and $\mu_j$, $j=1,2,...,m$. Since it is known that (Karlis, 2003, Larsson, 2020) $\hat\lambda+\hat\mu_j=\bar y_j$, the task is to maximize the likelihood in (\ref{L}) reformulated as
\begin{equation}
L(\lambda)=e^{-n\left\{\sum_{j=1}^m\bar y_j-(m-1)\lambda\right\}}\left\{\prod_{j=1}^m(\bar y_j-\lambda)^{n\bar y_j}\right\}\prod_{i=1}^n S_i(\lambda),\label{L1}
\end{equation}
where
\begin{equation}
S_i(\lambda)\eqd\sum_{u=0}^{y_i}\frac{\xi(\lambda)^u}{u!\prod_{j=1}^m(y_{ij}-u)!},\label{Si}
\end{equation}
with 
\begin{equation}
\xi(\lambda)\eqd\frac{\lambda}{\prod_{j=1}^m(\bar y_j-\lambda)}.\label{xi0}
\end{equation}

As usual, it seems to be easier to focus on maximizing the log likelihood, which via (\ref{L1}) is given as
\begin{align}
l(\lambda)&\eqd\log\{L(\lambda)\}\notag\\
&=n\left\{(m-1)\lambda-\sum_{j=1}^m\bar y_j\right\}+n\sum_{j=1}^m\bar y_j\log(\bar y_j-\lambda)+\sum_{i=1}^ns_i(\lambda),\label{logL}
\end{align}
where $s_i(\lambda)=\log\{S_i(\lambda)\}$. The derivative of (\ref{logL}) with respect to $\lambda$ is
\begin{equation}
l'(\lambda)=n\left(m-1-\sum_{j=1}^m\frac{\bar y_j}{\bar y_j-\lambda}\right)+\sum_{i=1}^ns_i'(\lambda).\label{lp}
\end{equation}

The derivative of $s_i(\lambda)$ w.r.t $\lambda$ is, emphasizing the dependency on the observations ${\boldsymbol y_i}=(y_{i1},y_{i2},...,y_{im})$,
$$s_i'(\lambda)=\frac{S_i'(\lambda;{\boldsymbol y_i})}{S_i(\lambda;{\boldsymbol y_i})},$$
where from (\ref{Si}), denoting $I\{A\}$ as the indicator of the event $A$,

\begin{align*}
&S_i'(\lambda;{\boldsymbol y_i})\\
&\eqd\xi'(\lambda)\sum_{u=1}^{y_i}\frac{\xi(\lambda)^{u-1}}{(u-1)!\prod_{j=1}^m(y_{ij}-u)!}I\{y_i>0\}\\
&=\xi'(\lambda)\sum_{u=1}^{(y_i-1)+1}\frac{\xi(\lambda)^{u-1}}{(u-1)!\prod_{j=1}^m\{(y_{ij}-1-(u-1)\}!}I\{y_i>0\}\\
&=\xi'(\lambda)S_i(\lambda;\boldsymbol{y}_i-\boldsymbol{1})I\{y_i>0\},
\end{align*}
where $\boldsymbol{1}=(1,1,...,1)$.
Hence, because by (\ref{xi0}),
$$\xi'(\lambda)=\frac{\prod_{j=1}^m(\bar y_j-\lambda)-\lambda\frac{d}{d\lambda}\prod_{j=1}^m(\bar y_j-\lambda)}{\left\{\prod_{j=1}^m(\bar y_j-\lambda)\right\}^2}
=\frac{1+\lambda\sum_{j=1}^m(\bar y_j-\lambda)^{-1}}{\prod_{j=1}^m(\bar y_j-\lambda)},
$$
it follows that
$$s_i'(\lambda)=\frac{1+\lambda\sum_{j=1}^m(\bar y_j-\lambda)^{-1}}{\prod_{j=1}^m(\bar y_j-\lambda)}\frac{S_i(\lambda;\boldsymbol{y}_i-\boldsymbol{1})}{S_i(\lambda;\boldsymbol{y_i})}I\{y_i>0\}.$$
Thus, via (\ref{lp}), and the identity
$$\sum_{j=1}^m\frac{\bar y_j}{\bar y_j-\lambda}=m+\lambda\sum_{j=1}^m(\bar y_j-\lambda)^{-1},$$
we have
\begin{align}
l'(\lambda)
&=\left\{1+\lambda\sum_{j=1}^m(\bar y_j-\lambda)^{-1}\right\}\notag\\
&\cdot\left\{-n+\frac{1}{\prod_{j=1}^m(\bar y_j-\lambda)}\sum_{i=1}^n
\frac{S_i(\lambda;\boldsymbol{y}_i-\boldsymbol{1})}{S_i(\lambda;\boldsymbol{y_i})}I\{y_i>0\}\right\}.\label{lp1}
\end{align}
To find the MLE $\hat\lambda$, we need to solve the equation $l'(\lambda)=0$. In practice, this has to be done by numerical methods. For our purposes, it is enough to derive an asymptotic approximation.

To this end, we Taylor expand (\ref{Si}) as (here, $y_{ij(2)}=y_{ij}(y_{ij}-1)$ etcetera, and $O$ in fact means order in probability)
\begin{align}
S_i(\lambda)=\frac{1}{\prod_{j=1}^m y_{ij}!}&\left\{1+\prod_{j=1}^my_{ij}\xi(\lambda)+\frac{\prod_{j=1}^m y_{ij(2)}}{2!}\xi(\lambda)^2+\frac{\prod_{j=1}^m y_{ij(3)}}{3!}\xi(\lambda)^3\right.\notag\\
&\left.+\frac{\prod_{j=1}^m y_{ij(4)}}{4!}\xi(\lambda)^4
+O(\lambda^5)\right\}.\label{Si2}
\end{align}
Here, Taylor expansion of (\ref{xi0}) yields (the sums and products are over $1\leq j\leq m$, etcetera), introducing $P=\prod_j\bar y_j$,
\begin{align}
&\xi(\lambda)\notag\\
&=\lambda P^{-1}\prod_{j=1}^m\left(1-\frac{\lambda}{\bar y_j}\right)^{-1}\notag\\
&=\lambda P^{-1}\prod_j\left(1+\frac{\lambda}{\bar y_j}+\frac{\lambda^2}{\bar y_j^2}+...\right)\notag\\
&=\lambda P^{-1}\left\{1+\sum_j\frac{1}{\bar y_j}\lambda+\left(\sum_j\frac{1}{\bar y_j^2}+\sum_{j<k}\frac{1}{\bar y_j\bar y_k}\right)\lambda^2\right.\notag\\
&\left.+\left(\sum_j\frac{1}{\bar y_j^3}+\sum_{j<k}\frac{1}{\bar y_j^2\bar y_k}+\sum_{j<k}\frac{1}{\bar y_j\bar y_k^2}+\sum_{j<k<l}\frac{1}{\bar y_j\bar y_k\bar y_l}\right)\lambda^3+O(\lambda^4)\right\},\label{xi2}
\end{align}
which inserted into (\ref{Si2}) implies
\begin{align*}
&S_i(\lambda)\\
&=\frac{1}{\prod_j y_{ij}!}\left[1+\prod_j y_{ij}\left\{P^{-1}\lambda+\frac{1}{2}a_2({\boldsymbol y}_i)P^{-2}\lambda^2+\frac{1}{6}a_3({\boldsymbol y}_i)P^{-3}\lambda^3
\right\}+O(\lambda^4)\right],
\end{align*}
where
\begin{equation}
a_2({\boldsymbol y}_i)\eqd2P\sum_j\bar y_j^{-1}+\prod_j (y_{ij}-1),\label{a2}
\end{equation}
and
\begin{align}
a_3({\boldsymbol y}_i)&\eqd 6P^2\left(\sum_j\bar y_j^{-2}+\sum_{j<k}\bar y_j^{-1}\bar y_k^{-1}\right)+6P\prod_j(y_{ij}-1)\sum_k\bar y_k^{-1}\notag\\
&+\prod_j(y_{ij}-1)(y_{ij}-2).\label{a3}
\end{align}
Now,
\begin{align}
&\frac{S_i(\lambda;\boldsymbol{y}_i-\boldsymbol{1})}{S_i(\lambda;\boldsymbol{y}_i)}=\prod_j y_{ij}\cdot\notag\\
&\frac{1+\prod_j(y_{ij}-1)\left\{P^{-1}\lambda+\frac{1}{2}a_2(\boldsymbol{y}_i-\boldsymbol{1})P^{-2}\lambda^2+
+\frac{1}{6}a_3(\boldsymbol{y}_i-\boldsymbol{1})P^{-3}\lambda^3\right\}+O(\lambda^4)}{1+\prod_j y_{ij}\left\{P^{-1}\lambda+\frac{1}{2}a_2(\boldsymbol{y}_i)P^{-2}\lambda^2+\frac{1}{6}a_3(\boldsymbol{y}_i)P^{-3}\lambda^3\right\}+O(\lambda^4)},\label{Sikvot}
\end{align}
and by Taylor expansion,
\begin{align}
\frac{S_i(\lambda;\boldsymbol{y}_i-\boldsymbol{1})}{S_i(\lambda;\boldsymbol{y}_i)}
&=\prod_j y_{ij}\left[1+\left\{\prod_j(y_{ij}-1)-\prod_j y_{ij}\right\}P^{-1}\lambda\right.\notag\\
&\left.+b_2(\boldsymbol{y}_i)P^{-2}\lambda^2+b_3(\boldsymbol{y}_i)P^{-3}\lambda^3
+O(\lambda^4)\right],\label{Sikvotexp}
\end{align}
where
\begin{align}
b_2(\boldsymbol{y}_i)
&\eqd\frac{1}{2}\left\{\prod_j(y_{ij}-1)a_2(\boldsymbol{y}_i-\boldsymbol{1})-\prod_j y_{ij}a_2(\boldsymbol{y}_i)\right\}\notag\\
&+\prod_j y_{ij}^2-\prod_j y_{ij}(y_{ij}-1).\label{b2}
\end{align}
and
\begin{align}
b_3(\boldsymbol{y}_i)&\eqd
\frac{1}{6}\left\{\prod_j(y_{ij}-1)a_3(\boldsymbol{y}_i-\boldsymbol{1})
-\prod_j y_{ij}a_3(\boldsymbol{y}_i)\right\}\notag\\
&-\prod_j y_{ij}^3+\prod_j y_{ij}^2a_2(\boldsymbol{y}_i)
+\prod_j y_{ij}^2(y_{ij}-1)\notag\\
&-\frac{1}{2}\prod_j y_{ij}\left\{\prod_j(y_{ij}-1)a_2(\boldsymbol{y}_i-\boldsymbol{1})
+\prod_j y_{ij}a_2(\boldsymbol{y}_i)\right\}.\label{b3}
\end{align}
Moreover, let
\begin{equation}
C_1\eqd n^{-1}\sum_{i=1}^n\prod_j y_{ij}\left\{\prod_j(y_{ij}-1)-\prod_j y_{ij}\right\},\label{C1}
\end{equation}
\begin{equation}
C_2\eqd n^{-1}\sum_{i=1}^n\prod_j y_{ij}b_2(\boldsymbol{y}_i),\label{C2}
\end{equation}
and
\begin{equation}
C_3\eqd n^{-1}\sum_{i=1}^n\prod_j y_{ij}b_3(\boldsymbol{y}_i),\label{C3}
\end{equation}
Now, via (\ref{Sikvotexp}) and (\ref{C1})-(\ref{C3}),
\begin{align}
S&\eqd n^{-1}\sum_{i=1}^n\frac{S_i(\lambda;\boldsymbol{y}_i-\boldsymbol{1})}{S_i(\lambda;\boldsymbol{y}_i)}I\{y_i>0\}\notag\\
&=\bar y_{11...1}+P^{-1}C_1\lambda
+P^{-2}C_2\lambda^2+P^{-3}C_3\lambda^3+O(\lambda^4),\label{S}
\end{align}
where
\begin{equation}
\bar y_{11...1}\eqd n^{-1}\prod_j y_{ij}.\label{y111}
\end{equation}

From (\ref{lp1}) and (\ref{S}), $l'(\lambda)=0$ implies
$\prod_j(\bar y_j-\lambda)=S$,
where
\begin{align*}
&\prod_j(\bar y_j-\lambda)=P\prod_j\left(1-\frac{\lambda}{\bar y_j}\right)\\
&=P\left(1-\sum_j\bar y_j^{-1}\lambda+\sum_{j<k}\bar y_j^{-1}\bar y_k^{-1}\lambda^2-\sum_{j<k<l}\bar y_j^{-1}\bar y_k^{-1}\bar y_l^{-1}\lambda^3\right)+O(\lambda^4),
\end{align*}
and so, via (\ref{S}),
\begin{equation}
0=n^{-1/2}d_0+d_1\lambda+d_2\lambda^2+d_3\lambda^3+O(\lambda^4),\label{lamsolve3}
\end{equation}
where
\begin{align}
d_0&=\sqrt n\left(\bar y_{11...1}-P\right),\label{d0}\\
d_1&=P^{-1}C_1+P\sum_j \bar y_j^{-1},\label{d1}\\
d_2&=P^{-2}C_2-P\sum_{j<k} \bar y_j^{-1}\bar y_k^{-1},\label{d2}\\
d_3&=P^{-3}C_3+P\sum_{j<k<l} \bar y_j^{-1}\bar y_k^{-1}\bar y_l^{-1}.\label{d3}
\end{align}
It will be shown below that as $n\to\infty$, $d_0$, $d_1$, $d_2$ and $d_3$ are $O_p(1)$. 

Now, putting $\lambda=n^{-1/2}c_1+n^{-1}c_2+n^{-3/2}c_3+O(n^{-2})$ and inserting into (\ref{lamsolve3}), we get
\begin{align*}
0&=n^{-1/2}d_0+d_1(n^{-1/2}c_1+n^{-1}c_2+n^{-3/2}c_3)+
d_2(n^{-1/2}c_1+n^{-1}c_2)^2\\
&+d_3(n^{-1/2}c_1)^3+O(n^{-2})\\
&=n^{-1/2}(d_0+d_1c_1)+n^{-1}(d_1c_2+d_2c_1^2)
+n^{-3/2}(d_1c_3+2d_2c_1c_2+d_3c_1^3)\\&+O(n^{-2}),
\end{align*}
and by putting all $n$ power coefficients equal to zero, we obtain the root
\begin{equation}
\hat\lambda=-\frac{n^{-1/2}d_0}{d_1}\left\{1+\frac{n^{-1/2}d_0d_2}{d_1^2}+\frac{n^{-1}d_0^2(2d_2^2-d_1d_3)}{d_1^4}+O(n^{-3/2})\right\}.\label{lam1}
\end{equation}

\section{The log likelihood ratio test}
Putting $\lambda=0$ in (\ref{L}) we get the likelihood under the null of independence. Then, inserting the MLEs $\hat\mu_j=\bar y_j$, we get the maximized likelihood under $H_0$ as
\begin{equation}
L_0\eqd e^{-n\sum_{j=1}^m\bar y_j}\left(\prod_{j=1}^m\bar y_j^{n\bar y_j}\right)\prod_{i=1}^n\prod_{j=1}^m (y_{ij}!)^{-1}.\label{L0}
\end{equation}
Now, define the log likelihood ratio test statistic as
$$Q_n=-2\log\left(\frac{L_0}{L_1}\right)=2(\log L_1-\log L_0),$$
where $L_1$ is the maximized likelihood under $H_1$. In the following, we will implicitly condition on the event $\hat\lambda>0$.

Via (\ref{logL}) and (\ref{L0}), we have
\begin{equation}
\frac{1}{2}n^{-1}Q_n=(m-1)\hat\lambda+\sum_j\bar y_j\log\left(1-\frac{\hat\lambda}{\bar y_j}\right)+n^{-1}\sum_{i=1}^n s_i^*(\hat\lambda),\label{Q}
\end{equation}
with
\begin{equation}
s_i^*(\hat\lambda)=\log\left\{S_i(\hat\lambda)\prod_j y_{ij}!\right\},\label{Sistar}
\end{equation}
where, via (\ref{Si2}),
\begin{align*}
&S_i(\hat\lambda)\prod_j y_{ij}!\\
&=1+\prod_j y_{ij}\xi(\hat\lambda)
+\frac{\prod_j y_{ij(2)}}{2!}\xi(\hat\lambda)^2
+\frac{\prod_j y_{ij(3)}}{3!}\xi(\hat\lambda)^3+\frac{\prod_j y_{ij(4)}}{4!}\xi(\hat\lambda)^4\\
&+O(\hat\lambda^5).
\end{align*}
Now, (\ref{Sistar}) and the Taylor expansion $\log(1+x)=x-x^2/2+x^3/3-x^4/4+O(x^5)$ yields
\begin{equation}
s_i^*(\hat\lambda)=\prod_j y_{ij}\xi(\hat\lambda)+\frac{e_{2i}}{2}\xi(\hat\lambda)^2+\frac{e_{3i}}{3}\xi(\hat\lambda)^3+\frac{e_{4i}}{4}\xi(\hat\lambda)^4+O(\hat\lambda^5),
\label{sistarexp}
\end{equation}
where
\begin{align}
e_{2i}&=\prod_j y_{ij(2)}-\left(\prod_j y_{ij}\right)^2,\label{e2i}\\
e_{3i}&=\frac{1}{2}\prod_j y_{ij(3)}-\frac{3}{2}\prod_j y_{ij}\prod_j y_{ij(2)}
+\left(\prod_j y_{ij}\right)^3,\label{e3i}\\
e_{4i}&=\frac{1}{6}\prod_j y_{ij(4)}-\frac{2}{3}\prod_j y_{ij}\prod_j y_{ij(3)}
-\frac{1}{2}\left(\prod_j y_{ij(2)}\right)^2\notag\\
&+2\left(\prod_j y_{ij}\right)^2\prod_j y_{ij(2)}-\left(\prod_j y_{ij}\right)^4.\label{e4i}
\end{align}
Moreover, because for $j=1,2,...,m$,
$$\log\left(\frac{\bar y_j-\hat\lambda}{\bar y_j}\right)=\log\left(1-\frac{\hat\lambda}{\bar y_j}\right)
=-\frac{\hat\lambda}{\bar y_j}-\frac{\hat\lambda^2}{2\bar y_j^2}-\frac{\hat\lambda^3}{3\bar y_j^3}-\frac{\hat\lambda^4}{4\bar y_j^4}+O(\hat\lambda^5),
$$
(\ref{Q}) and (\ref{sistarexp}) yield
\begin{align}
&\frac{1}{2}n^{-1}Q_n\notag\\&=(m-1)\hat\lambda+\sum_{j=1}^m\bar y_j\left(-\frac{\hat\lambda}{\bar y_j}-\frac{\hat\lambda^2}{2\bar y_j^2}-\frac{\hat\lambda^3}{3\bar y_j^3}-\frac{\hat\lambda^4}{4\bar y_j^4}\right)+\bar y_{11...1}\xi(\hat\lambda)\notag\\
&+\frac{1}{2}\bar e_2\xi(\hat\lambda)^2+\frac{1}{3}\bar e_3\xi(\hat\lambda)^3+\frac{1}{4}\bar e_4\xi(\hat\lambda)^4+O(\hat\lambda^5),
\label{Qb}
\end{align}
where
$\bar e_k=n^{-1}\sum_i e_{ki}$ for $k=2,3,4$.
Here, from (\ref{e2i}) and (\ref{C1}),
\begin{equation}
\bar e_2=C_1,\label{e2bar}\end{equation}

Let $\by(2^{m-1},0^1)$ be the sum of all products of $m-1$ distinct ${\bar y_j}^2$, etcetera.
Then, inserting (\ref{xi2}) into (\ref{Qb}), we get in view of (\ref{d0}) that
\begin{align*}
&\frac{1}{2}n^{-1}Q_n\notag\\
&=-\hat\lambda-\frac{\by(1^{m-1},0^1)}{2\prod_j\bar y_j}\hat\lambda^2-\frac{\by(2^{m-1},0^1)}{3\prod_j\bar y_j^2}\hat\lambda^3-\frac{\by(3^{m-1},0^1)}{4\prod_j\bar y_j^3}\hat\lambda^4\\
&+\bar y_{11...1}\left(\frac{1}{\prod_j\bar y_j}\hat\lambda+\frac{\by(1^{m-1},0^1)}{\prod_j\bar y_j^2}\hat\lambda^2+\frac{\by(2^{m-1},0^1)+\by(2^{m-2},1^2)}{\prod_j\bar y_j^3}\hat\lambda^3\right.\\
&\left.+\frac{\by(3^{m-1},0^1)+\by(3^{m-2},2^1,1^1)+\by(3^{m-3},2^3)}{\prod_j\bar y_j^4}\hat\lambda^4\right)\\
&+\frac{1}{2}\bar e_2\left(\frac{1}{\prod_j\bar y_j}\hat\lambda+\frac{\by(1^{m-1},0^1)}{\prod_j\bar y_j^2}\hat\lambda^2+\frac{\by(2^{m-1},0^1)+\by(2^{m-2},1^2)}{\prod_j\bar y_j^3}\hat\lambda^3\right)^2\\
&+\frac{1}{3}\bar e_3\left(\frac{1}{\prod_j\bar y_j}\hat\lambda+\frac{\by(1^{m-1},0^1)}{\prod_j\bar y_j^2}\hat\lambda^2\right)^3+\frac{1}{4}\bar e_4\left(\frac{1}{\prod_j\bar y_j}\hat\lambda\right)^4+O(\hat\lambda^5)\\
&=\frac{n^{-1/2}d_0}{\prod_j\bar y_j}\hat\lambda
+\frac{d_{32}}{2\prod_j\bar y_j^2}\hat\lambda^2
+\frac{d_{33}}{3\prod_j\bar y_j^3}\hat\lambda^3+\frac{d_{34}}{4\prod_j\bar y_j^4}\hat\lambda^4+O(\hat\lambda^5),
\end{align*}
where by some simplifications,
\begin{align}
d_{32}&=(2\bar y_{11...1}-P)\by(1^{m-1},0^1)+\bar e_2,\label{d32}\\
d_{33}&=(3\bar y_{11...1}-P)\by(2^{m-1},0^1)
+3\bar y_{11...1}\by(2^{m-2},1^2)\notag\\
&+3\by(1^{m-1},0^1)\bar e_2+\bar e_3,
\label{d33}\\
d_{34}&=(4\bar y_{11...1}-P)\by(3^{m-1},0^1)
\notag\\
&+4\bar y_{11...1}\left\{\by(3^{m-2},2^1,1^1)+\by(3^{m-3},2^3)\right\}\notag\\
&+2\left\{\by(1^{m-1},0^1)^2+2\by(2^{m-1},0^1)+2\by(2^{m-2},1^2)\right\}\bar e_2\notag\\
&+4\by(1^{m-1},0^1)\bar e_3+\bar e_4.\label{d34}
\end{align}
Hence, by (\ref{lam1}),
\begin{align}
\frac{1}{2}Q
&=\frac{d_0^2}{P(-d_1)}\left\{1+\frac{n^{-1/2}d_0d_2}{d_1^2}+\frac{n^{-1}d_0^2(2d_2^2-d_1d_3)}{d_1^4}\right\}\notag\\
&+\frac{d_{32}d_0^2}{2P^2d_1^2}\left\{1+\frac{n^{-1/2}d_0d_2}{d_1^2}+\frac{n^{-1}d_0^2(2d_2^2-d_1d_3)}{d_1^4}\right\}^2\notag\\
&+\frac{n^{-1/2}d_{33}d_0^3}{3P^3(-d_1)^3}\left(1+\frac{n^{-1/2}d_0d_2}{d_1^2}\right)^3
+\frac{n^{-1}d_{34}d_0^4}{4P^4d_1^4}
+O(n^{-3/2})\notag\\
&=f_0+n^{-1/2}f_1+n^{-1}f_2+O(n^{-3/2}),\label{Qc}
\end{align}
where
\begin{align}
f_0&=\frac{d_0^2}{2P^2d_1^2}(-2Pd_1+d_{32}),\label{f0}\\
f_1&=\frac{d_0^3}{3P^3d_1^4}(-3P^2d_1d_2+3Pd_2d_{32}-d_{33}d_1),\label{f1}\\
f_2&=\frac{d_0^4}{4P^4d_1^6}(-8P^3d_1d_2^2+10P^2d_2^2d_{32}-4Pd_{33}d_1d_2+d_{34}d_1^2\notag\\
&\qquad\qquad+4P^3d_1^2d_3-4P^2d_1d_3d_{32}).\label{f2}
\end{align}

In order to proceed from here, for $j=1,2,...,m$ write 
\begin{equation}
\bar y_j=\mu_j+n^{-1/2}r_j,\label{baryj}
\end{equation}
where
\begin{equation}
r_j=\sqrt n(\bar y_j-\mu_j).\label{rj}
\end{equation}
At first, we  use (\ref{baryj}) and the expansion $(1+x)^{-1}=1-x+x^2+O(x^3)$ to get
\begin{align}
\bar y_j^{-1}&=(\mu_j+n^{-1/2}r_j)^{-1}=\mu_j^{-1}\left(1+n^{-1/2}\frac{r_j}{\mu_j}\right)^{-1}\notag\\
&=\mu_j^{-1}\left(1-n^{-1/2}\frac{r_j}{\mu_j}+n^{-1}\frac{r_j^2}{\mu_j^2}\right)+O(n^{-3/2}).\label{baryjexp}
\end{align}
Now, defining $p_1=\prod_j\mu_j$, we have because $E(P)=p_1$ that, via (\ref{baryj}),
\begin{align}
P&=p_1\prod_j\left(1+n^{-1/2}\frac{r_j}{\mu_j}\right)\notag\\
&=p_1\left\{1+n^{-1/2}\sum_j\frac{r_j}{\mu_j}+n^{-1}\sum_{j<k}\frac{r_jr_k}{\mu_j\mu_k}+O(n^{-3/2})\right\},\label{P}
\end{align}
implying
\begin{align}
&P^{-1}\notag\\
&=p_1^{-1}\left[1-n^{-1/2}\sum_j\frac{r_j}{\mu_j}+n^{-1}\left\{\left(\sum_j\frac{r_j}{\mu_j}\right)^2-\sum_{j<k}\frac{r_jr_k}{\mu_j\mu_k}\right\}+O(n^{-3/2})\right]\notag\\
&=p_1^{-1}\left\{1-n^{-1/2}\sum_j\frac{r_j}{\mu_j}+n^{-1}\sum_{j\geq k}\frac{r_jr_k}{\mu_j\mu_k}+O(n^{-3/2})\right\}.
\label{Pinv}
\end{align}
To tackle $d_1=P^{-1}C_1+P\sum_j\bar y_j^{-1}$ (cf (\ref{d1})), we write 
\begin{equation}
C_1=C_{10}+n^{-1/2}C_{11},\quad C_{11}=\sqrt{n}(C_1-C_{10}),\label{C1b}
\end{equation}
where $C_{10}=E(C_1)$. 
Now, write 
\begin{equation}
-d_1=d_{10}+n^{-1/2}d_{11}+n^{-1}d_{12}+O(n^{-3/2}),\label{md1}
\end{equation}
where from (\ref{d1}), for $m>2$,
\begin{equation}
-d_{10}=p_1^{-1}C_{10}+p_1s_1,\label{d10}
\end{equation}
letting $s_1\eqd\sum_j\mu_j^{-1}$.

Now, from (\ref{d1}) and (\ref{baryjexp})-(\ref{C1b}),
\begin{align*}
d_1&=p_1^{-1}\left(1-n^{-1/2}\sum_j\frac{r_j}{\mu_j}+n^{-1}\sum_{j\geq k}\frac{r_jr_k}{\mu_j\mu_k}\right)(C_{10}+n^{-1/2}C_{11})\\
&+p_1\left(1+n^{-1/2}\sum_j\frac{r_j}{\mu_j}+n^{-1}\sum_{j< k}\frac{r_jr_k}{\mu_j\mu_k}\right)\sum_j\mu_j^{-1}\left(1-n^{-1/2}\frac{r_j}{\mu_j}+n^{-1}\frac{r_j^2}{\mu_j^2}\right)\\
&+O(n^{-3/2}),
\end{align*}
from which it follows that
\begin{align}
d_{11}&=p_1^{-1}\left(C_{10}\sum_j\frac{r_j}{\mu_j}-C_{11}\right)-p_1\left(\sum_j\mu_j^{-1}\sum_j\frac{r_j}{\mu_j}-\sum_j\frac{r_j}{\mu_j^2}\right),\label{d11}\\
d_{12}&=p_1^{-1}\left(C_{11}\sum_j\frac{r_j}{\mu_j}-C_{10}\sum_{j\geq k}\frac{r_jr_k}{\mu_j\mu_k}\right)\notag\\
&-p_1\left(\sum_j\mu_j^{-1}\sum_{j< k}\frac{r_jr_k}{\mu_j\mu_k}-\sum_j\frac{r_j}{\mu_j}\sum_j\frac{r_j}{\mu_j^2}+\sum_j\frac{r_j^2}{\mu_j^3}\right).\label{d12}
\end{align}
Moreover, from (\ref{md1}),
\begin{equation}
(-d_1)^{-1}=d_{10}^{-1}\left(1-n^{-1/2}\frac{d_{11}}{d_{10}}+n^{-1}\frac{d_{11}^2-d_{10}d_{12}}{d_{10}^2}\right)+O(n^{-3/2}).\label{md1inv}
\end{equation}

Now, write (\ref{f0}) as
\begin{equation}
f_0=f_{01}+f_{02},\label{f0dec}
\end{equation}
where
\begin{align}
f_{01}&=\frac{d_0^2}{P(-d_1)},\label{f01}\\
f_{02}&=\frac{d_0^2d_{32}}{2P^2d_1^2}.\label{f02}
\end{align}
Next, introduce the notation $\dot y_{ij}=y_{ij}-\mu_j$ for $j=1,2,...m$. Moreover, let $\bar{\dot y}_j=n^{-1}\sum_i\dot y_{ij}$, $\bar{\dot y}_{110...0}=n^{-1}\sum_i\dot y_{1j}y_{2j}$ and so on. It follows from (\ref{d0}) that, for $m>2$,
\begin{align}
d_0&=\sqrt{n}(\bar y_{11...1}-\bar y_1\bar y_2\cdots\bar y_m)\notag\\
&=n^{-1/2}\sum_i(\mu_1+\dot y_{i1})\cdots(\mu_m+\dot y_{im})-n^{1/2}(\mu_1+\bar{\dot y}_1)\cdots(\mu_m+\bar{\dot y}_m)\notag\\
&=\sqrt{n}(\bar{\dot y}_{11...1}+\mu_1\bar{\dot y}_{01...1}+...+\mu_m\bar{\dot y}_{1...10}
+...\notag\\
&+\mu_1\mu_2\cdots\mu_{m-2}\bar{\dot y}_{0...011}+...+\mu_3\cdots\mu_m\bar{\dot y}_{110...0})\notag\\
&-\sqrt{n}(\mu_1\mu_2\cdots\mu_{m-2}\bar{\dot y}_{m-1}\bar{\dot y}_m+...+\mu_3\cdots\mu_m\bar{\dot y}_1\bar{\dot y}_2)\notag\\
&-\sqrt{n}(\mu_1\mu_2\cdots\mu_{m-3}\bar{\dot y}_{m-2}\bar{\dot y}_{m-1}\bar{\dot y}_m+...+\mu_4\cdots\mu_m\bar{\dot y}_1\bar{\dot y}_2\bar{\dot y}_3)+...\notag\\
&=d_{00}-n^{-1/2}(r^2,\mu^{m-2})-n^{-1}(r^3,\mu^{m-3})+O(n^{-3/2}),\label{d0b}
\end{align} 
where for $k=2,3$ (note that from (\ref{rj}), $r_j=\sqrt{n}\bar{\dot y}_j$)
\begin{equation}
(r^k,\mu^{m-k})=\mu_1\mu_2\cdots\mu_{m-k}r_{m-k+1}\cdots r_m+...+\mu_{k+1}\cdots\mu_m r_1\cdots r_k,\label{rmu}
\end{equation}
and
\begin{equation}
d_{00}=\sqrt{n}\bar{\dot y}_{11...1}+(r^{\otimes(m-1)},\mu^1)+...+(r^{\otimes 2},\mu^{m-2}),\label{d00}
\end{equation}
with, for $k=1,2,...,m-2$,
\begin{equation}
(r^{\otimes(m-k)},\mu^k)=\sqrt{n}(\mu_1\mu_2\cdots\mu_k\bar{\dot y}_{0...01...1}+...+\mu_{m-k+1}\cdots\mu_m\bar{\dot y}_{1...10...0}),
\label{rmu2}
\end{equation}
where the $\bar{\dot y}$ terms have indices with $k$ zeros.
It is seen that $d_{00}$ and the $(r^k,\mu^{m-k})$ are $O(1)$.

If $m=2$, it follows that $d_0=d_{00}-n^{-1/2}r_1r_2$, where $d_{00}=\sqrt{n}\bar{\dot y}_{11}$.

Furthermore, (\ref{f01}), (\ref{Pinv}) and (\ref{md1inv}) yield
\begin{align}
f_{01}
&=\frac{d_0^2}{p_1d_{10}}\left(1-n^{-1/2}\sum_j\frac{r_j}{\mu_j}+n^{-1}\sum_{j\geq k}\frac{r_jr_k}{\mu_j\mu_k}\right)\notag\\
&\left(1-n^{-1/2}\frac{d_{11}}{d_{10}}+n^{-1}\frac{d_{11}^2-d_{10}d_{12}}{d_{10}^2}\right)+O(n^{-3/2})\notag\\
&=\frac{d_0^2}{p_1d_{10}}(1+n^{-1/2}s_{011}+n^{-1}s_{012}),\label{f01b}
\end{align}
where
\begin{align}
s_{011}&=-\left(\sum_j\frac{r_j}{\mu_j}+\frac{d_{11}}{d_{10}}\right),\label{s011}\\
s_{012}&=\sum_{j\geq k}\frac{r_jr_k}{\mu_j\mu_k}+\frac{d_{11}}{d_{10}}\sum_j\frac{r_j}{\mu_j}
+\frac{d_{11}^2-d_{10}d_{12}}{d_{10}^2}.\label{s012}
\end{align}
To calculate the expectation of (\ref{f01b}), we at first have, via (\ref{d0b}),
\begin{align}
E(d_0^2)&=E(d_{00}^2)-2n^{-1/2}E\{d_{00}(r^2,\mu^{m-2})\}\notag\\
&+n^{-1}\left[E\{(r^2,\mu^{m-2})^2\}-2E\{d_{00}(r^3,\mu^{m-3})\}\right]+O(n^{-3/2}),\label{Ed02}
\end{align}
where, from (\ref{d00}),
\begin{align}
&E(d_{00}^2)\notag\\
&=nE(\bar{\dot y}_{11...1}^2)+E\left\{(r^{\otimes(m-1)},\mu^1)^2\right\}+...+E\left\{(r^{\otimes 2},\mu^{m-2})^2\right\}\notag\\
&+2\sqrt{n}E\left\{\bar{\dot y}_{11...1}(r^{\otimes(m-1)},\mu^1)\right\}
+...+E\left\{(r^{\otimes 3},\mu^{m-3})(r^{\otimes 2},\mu^{m-2})\right\}.
\label{Ed002}
\end{align}
Here, because all $E(\dot y_{ij})=0$,
\begin{align}
nE(\bar{\dot y}_{11...1}^2)&=n^{-1}E\left\{\left(\sum_i \dot y_{i1}\cdots \dot y_{im}\right)^2\right\}\notag\\
&=n^{-1}\sum_i\sum_j E(\dot y_{i1}\cdots \dot y_{im}\dot y_{j1}\cdots \dot y_{jm})\notag\\
&=n^{-1}\sum_i E(\dot y_{i1}^2)\cdots E(\dot y_{im}^2)=\mu_1\cdots\mu_m=p_1.\label{Ebary2}
\end{align}
A similar calculation shows that $E(\bar{\dot y}_I\bar{\dot y}_J)=0$ for index sets $I\neq J$. As a consequence, we have for $k=1,...,m-2$ that, from (\ref{rmu2}),
\begin{align*}
&E\left\{(r^{\otimes(m-k)},\mu^k)^2\right\}\\
&=nE\left\{(\mu_1\mu_2\cdots\mu_k\bar{\dot y}_{0...01...1}+...+\mu_{m-k+1}\cdots\mu_m\bar{\dot y}_{1...10...0})^2\right\}\\
&=n\left\{\mu_1^2\cdots \mu_k^2E(\bar{\dot y}_{0...01...1}^2)
+...+\mu_{m-k+1}^2\cdots \mu_m^2E(\bar{\dot y}_{1...10...0}^2)\right\},
\end{align*}
where as in (\ref{Ebary2}),
$$nE(\bar{\dot y}_{1...10...0}^2)=\mu_1\cdots\mu_{m-k}$$
etcetera, and it follows that
\begin{equation}
E\left\{(r^{\otimes(m-k)},\mu^k)^2\right\}=p_1\sum_{j_1<...<j_k}\mu_{j_1}\cdots\mu_{j_k}.\label{Ermu2}
\end{equation}
Moreover, as above,
$$\sqrt{n}E\left\{\bar{\dot y}_{11...1}(r^{\otimes(m-1)},\mu^1)\right\}
=nE\left\{\bar{\dot y}_{11...1}(\mu_1\bar{\dot y}_{01...1}
+...+\mu_m\bar{\dot y}_{1...10})\right\}=0,
$$
etcetera, and so, via (\ref{Ed002})-(\ref{Ermu2}) and (\ref{d10}), we find that
\begin{equation}
E(d_{00}^2)=p_1d_{10}.\label{Ed002b}
\end{equation}
As for the next term in (\ref{Ed02}), we have
\begin{align}
&n^{-1/2}E\{d_{00}(r^2,\mu^{m-2})\}\notag\\
&=n^{-1/2}E\left[\left\{\sqrt{n}\bar{\dot y}_{11...1}+(r^{\otimes m-1},\mu^1)+...+(r^{\otimes 2},\mu^{m-2})\right\}(r^2,\mu^{m-2})\right],\label{Ed022}
\end{align}
where
\begin{align*}
E\left\{\bar{\dot y}_{11...1}(r^2,\mu^{m-2})\right\}
&=E\left\{\bar{\dot y}_{11...1}(\mu_1\cdots\mu_{m-2}r_{m-1}r_m+...
+\mu_3\cdots\mu_mr_1r_2)\right\}\\
&=\mu_1\cdots\mu_{m-2}E(\bar{\dot y}_{11...1}r_{m-1}r_m)+...
+\mu_3\cdots\mu_mE(\bar{\dot y}_{11...1}r_1r_2).
\end{align*}
Here,
$$E(\bar{\dot y}_{11...1}r_1r_2)=n^{-2}\sum_i\sum_j\sum_k E(\dot y_{i1}\cdots\dot y_{im}\dot y_{j1}\dot y_{k2}),$$
If $m>2$ this is zero. For $m=2$, it is 
$$n^{-2}\sum_iE(\dot y_{i1}^2\dot y_{i2}^2)=n^{-1}\mu_1\mu_2.$$
In this case, this is the only contribution to (\ref{Ed022}). For $m>2$, we also have terms like
\begin{align*}
E\left\{(r^{\otimes 2},\mu^{m-2})(r^2,\mu^{m-2})\right\}=n^{3/2}E&\left\{(\mu_1\cdots\mu_{m-2}\bar{\dot y}_{0...011}+...+\mu_3\cdots\mu_m\bar{\dot y}_{110...0})\right.\\
&\left.(\mu_1\cdots\mu_{m-2}\bar{\dot y}_{m-1}\bar{\dot y}_m+...+\mu_3\cdots\mu_m\bar{\dot y}_1\bar{\dot y}_2)\right\},
\end{align*}
where to evaluate the right-hand side, we need terms like
\begin{align*}
E(\bar{\dot y}_{110...0}\bar{\dot y}_1\bar{\dot y}_2)
&=n^{-3}\sum_i\sum_j\sum_k E(\dot y_{i1}\dot y_{i2}\dot y_{j1}\dot y_{k2})=n^{-2}\mu_1\mu_2,\\
E(\bar{\dot y}_{110...0}\bar{\dot y}_1\bar{\dot y}_3)
&=n^{-3}\sum_i\sum_j\sum_k E(\dot y_{i1}\dot y_{i2}\dot y_{j1}\dot y_{k3})=0,
\end{align*}
etcetera. It follows that
\begin{align}
&E\left\{(r^{\otimes 2},\mu^{m-2})(r^2,\mu^{m-2})\right\}\notag\\
&=n^{-1/2}(\mu_1^2\cdots\mu_{m-2}^2\mu_{m-1}\mu_m+...
+\mu_3^2\cdots\mu_m^2\mu_1\mu_2)\notag\\
&=n^{-1/2}p_1\sum_{j_1<...<j_{m-2}}\mu_{j_1}\cdots\mu_{j_{m-2}}.\label{Ert2mr2m}
\end{align}
Similarly,
\begin{align*}
E\left\{(r^{\otimes 3},\mu^{m-3})(r^2,\mu^{m-2})\right\}=n^{3/2}E&\left\{(\mu_1\cdots\mu_{m-3}\bar{\dot y}_{0...0111}+...+\mu_4\cdots\mu_m\bar{\dot y}_{1110...0})\right.\\
&\left.(\mu_1\cdots\mu_{m-2}\bar{\dot y}_{m-1}\bar{\dot y}_m+...+\mu_3\cdots\mu_m\bar{\dot y}_1\bar{\dot y}_2)\right\},
\end{align*}
where
$$E(\bar{\dot y}_{1110...0}\bar{\dot y}_1\bar{\dot y}_2)
=n^{-3}\sum_i\sum_j\sum_k E(\dot y_{i1}\dot y_{i2}\dot y_{i3}\dot y_{j1}\dot y_{k2})=0,
$$
etcetera, implying 
$$E\left\{(r^{\otimes 3},\mu^{m-3})(r^2,\mu^{m-2})\right\}=0.$$
In this manner, (\ref{Ed022}) and (\ref{Ert2mr2m}) imply
\begin{align}
n^{-1/2}E\left\{d_{00}(r^2,\mu^{m-2})\right\}
&=n^{-1/2}E\left\{(r^{\otimes 2},\mu^{m-2})(r^2,\mu^{m-2})\right\}
\notag\\
&=n^{-1}p_1\sum_{j_1<...<j_{m-2}}\mu_{j_1}\cdots\mu_{j_{m-2}}.\label{Ed023}
\end{align}
To tackle the $n^{-1}$ terms of (\ref{Ed02}), we have
$$E\left\{(r^2,\mu^{m-2})^2\right\}
=n^2E\left\{(\mu_1\cdots\mu_{m-2}\bar{\dot y}_{m-1}\bar{\dot y}_m+...+\mu_3\cdots\mu_m\bar{\dot y}_1\bar{\dot y}_2)^2\right\},
$$
where
\begin{align*}
E(\bar{\dot y}_1^2\bar{\dot y}_2^2)
&=n^{-4}\sum_i\sum_j\sum_k\sum_lE(\dot y_{i1}\dot y_{j1}\dot y_{k2}\dot y_{l2})=n^{-2}\mu_1\mu_2+O(n^{-3}),\\
E(\bar{\dot y}_1^2\bar{\dot y}_2\bar{\dot y}_3)
&=n^{-4}\sum_i\sum_j\sum_k\sum_lE(\dot y_{i1}\dot y_{j1}\dot y_{k2}\dot y_{l3})=0,
\end{align*}
and so on, implying
\begin{align}
E\left\{(r^2,\mu^{m-2})^2\right\}&=\mu_1^2\cdots\mu_{m-2}^2\mu_{m-1}\mu_m+...+\mu_3^2\cdots\mu_m^2\mu_1\mu_2\notag\\
&=p_1\sum_{j_1<...<j_{m-2}}\mu_{j_1}\cdots\mu_{j_{m-2}}.\label{Ed031}
\end{align}
Moreover, from (\ref{d00}),
\begin{align*}
&E\{d_{00}(r^3,\mu^{m-3})\}\notag\\
&=E\left[\left\{\sqrt{n}\bar{\dot y}_{11...1}+(r^{\otimes m-1},\mu^1)+...+(r^{\otimes 2},\mu^{m-2})\right\}(r^3,\mu^{m-3})\right],
\end{align*}
where as above,
\begin{align*}
&E\left\{(r^{\otimes 2},\mu^{m-2})(r^3,\mu^{m-3})\right\}\\
&=n^2E\left\{(\mu_1\cdots\mu_{m-2}\bar{\dot y}_{0...011}+...+\mu_3\cdots\mu_m\bar{\dot y}_{110...0})\right.\\
&\left.(\mu_1\cdots\mu_{m-3}\bar{\dot y}_{m-2}\bar{\dot y}_{m-1}\bar{\dot y}_m+...+\mu_4\cdots\mu_m\bar{\dot y}_1\bar{\dot y}_2\bar{\dot y}_3)\right\}\\
&=0,
\end{align*}
and so on. The only term that breaks this pattern is
\begin{align*}
&E\left\{(r^{\otimes 3},\mu^{m-3})(r^3,\mu^{m-3})\right\}\\
&=n^2E\left\{(\mu_1\cdots\mu_{m-3}\bar{\dot y}_{0...0111}+...+\mu_4\cdots\mu_m\bar{\dot y}_{1110...0})\right.\\
&\left.(\mu_1\cdots\mu_{m-3}\bar{\dot y}_{m-2}\bar{\dot y}_{m-1}\bar{\dot y}_m+...+\mu_4\cdots\mu_m\bar{\dot y}_1\bar{\dot y}_2\bar{\dot y}_3)\right\},
\end{align*}
where nonzero contributions arise from terms like
$$E(\bar{\dot y}_{1110...0}\bar{\dot y}_1\bar{\dot y}_2\bar{\dot y}_3)
=n^{-4}\sum_i\sum_j\sum_k\sum_l E({\dot y}_{i1}{\dot y}_{i2}{\dot y}_{i3}{\dot y}_{j1}{\dot y}_{k2}{\dot y}_{l3})
=n^{-3}\mu_1\mu_2\mu_3.
$$
But this means that
$$E\{d_{00}(r^3,\mu^{m-3})\}=O(n^{-1}),$$
which may be neglected.  Thus, (\ref{Ed02}), (\ref{Ed002b}), (\ref{Ed023}) and (\ref{Ed031}) imply
\begin{equation}
E(d_0^2)=p_1\left(d_{10}-n^{-1}\sum_{j_1<...<j_{m-2}}\mu_{j_1}\cdots\mu_{j_{m-2}}\right)+O(n^{-3/2}),\label{Ed02b}
\end{equation}
for $m>2$, and the equation is true also for $m=2$ if the sum is interpreted as 1.

Next, we go on with $E(d_0^2 s_{011})$. To this end, using (\ref{d0b}) (cf (\ref{Ed02})),
\begin{align}
E(d_0^2r_1)&=E(d_{00}^2r_1)-2n^{-1/2}E\{d_{00}(r^2,\mu^{m-2})r_1\}\notag\\
&+n^{-1}\left[E\{(r^2,\mu^{m-2})^2r_1\}-2E\{d_{00}(r^3,\mu^{m-3})r_1\}\right]+O(n^{-3/2}),\label{Ed02r1}
\end{align}
where, from (\ref{d00}) (cf (\ref{Ed002})),
\begin{align}
&E(d_{00}^2r_1)\notag\\
&=nE(\bar{\dot y}_{11...1}^2r_1)+E\left\{(r^{\otimes(m-1)},\mu^1)^2r_1\right\}+...+E\left\{(r^{\otimes 2},\mu^{m-2})^2r_1\right\}\notag\\
&+2\sqrt{n}E\left\{\bar{\dot y}_{11...1}(r^{\otimes(m-1)},\mu^1)r_1\right\}
+...+E\left\{(r^{\otimes 3},\mu^{m-3})(r^{\otimes 2},\mu^{m-2})r_1\right\}.
\label{Ed002r1}
\end{align}
Here, in the usual manner,
\begin{align}
nE(\bar{\dot y}_{11...1}^2r_1)
&=n^{3/2}E(\bar{\dot y}_{11...1}^2\bar{\dot y}_1)\notag\\
&=n^{-3/2}\sum_i\sum_j\sum_k E(\dot y_{i1}\cdots\dot y_{im}\dot y_{j1}\cdots\dot y_{jm}\dot y_{k1})\notag\\
&=n^{-3/2}\sum_iE(\dot y_{i1}^3)E(\dot y_{i2}^2)\cdots E(\dot y_{im}^2)\notag\\
&=n^{-1/2}p_1,\label{Ebary2r1}
\end{align}
since the third central moment of a ${\rm Po}(\mu)$ variate is $\mu$. Moreover, we have for $k=1,...,m-2$ that
\begin{align*}
&E\left\{(r^{\otimes(m-k)},\mu^k)^2r_1\right\}\\
&=n^{3/2}E\left\{(\mu_1\mu_2\cdots\mu_k\bar{\dot y}_{0...01...1}+...+\mu_{m-k+1}\cdots\mu_m\bar{\dot y}_{1...10...0})^2\bar{\dot y}_1\right\}\\
&=n^{3/2}\left\{\mu_1^2\cdots \mu_k^2E(\bar{\dot y}_{0...01...1}^2\bar{\dot y}_1)
+...+\mu_{m-k+1}^2\cdots \mu_m^2E(\bar{\dot y}_{1...10...0}^2\bar{\dot y}_1)\right\},
\end{align*}
because as above, cross terms disappear in the last step. Here, the expectations are nonzero only if the first index in $\bar{\dot y}_I$ is one, and in such cases we e.g. have
\begin{align}
E(\bar{\dot y}_{1...10...0}^2\bar{\dot y}_1)
&=n^{-3}\sum_i\sum_j\sum_k E(\dot y_{i1}\cdots\dot y_{i,m-k}\dot y_{j1}\cdots\dot y_{j,m-k}\dot y_{k1})\notag\\
&=n^{-2}\mu_1\cdots\mu_{m-k}.\label{Ebary2bary1}
\end{align}
Hence,
\begin{align}
&E\left\{(r^{\otimes(m-k)},\mu^k)^2r_1\right\}\notag\\
&=n^{3/2}\left\{\mu_2^2\cdots\mu_{k+1}^2E(\bar{\dot y}_{10...01...1}^2\bar{\dot y}_1)+...+\mu_{m-k+1}^2\cdots\mu_m^2E(\bar{\dot y}_{1...10...0}^2\bar{\dot y}_1)\right\}\notag\\
&=n^{-1/2}\left(\mu_2^2\cdots\mu_{k+1}^2\mu_1\mu_{k+1}\cdots\mu_m+...+\mu_{m-k+1}^2\cdots\mu_m^2\mu_1\cdots\mu_{m-k}\right)
\notag\\
&=n^{-1/2}p_1\sum_{1<j_1<...<j_k}\mu_{j_1}\cdots\mu_{j_k}.
\label{Ertmkmkr1}
\end{align}
For the usual reasons, there is only one more term that contributes to (\ref{Ed002r1}), and this term is
\begin{align*}
&E\left\{(r^{\otimes 3},\mu^{m-3})(r^{\otimes 2},\mu^{m-2})r_1\right\}
\\
&=n^{3/2}E\left\{(\mu_1\cdots\mu_{m-3}\bar{\dot y}_{0...0111}+...
+\mu_4\cdots\mu_m\bar{\dot y}_{1110...0})\right.\\
&\left.(\mu_1\cdots\mu_{m-2}\bar{\dot y}_{0...011}+...
+\mu_3\cdots\mu_m\bar{\dot y}_{110...0})\bar{\dot y}_1\right\}.
\end{align*}
The nonzero expectations that arise from here are of the type
\begin{align*}
E(\bar{\dot y}_{1110...0}\bar{\dot y}_{0110...0}\bar{\dot y}_1)
&=n^{-3}\sum_i\sum_j\sum_k E(\dot y_{i1}\dot y_{i2}\dot y_{i3}\dot y_{j2}\dot y_{j3}\dot y_{k1})\\
&=n^{-3}\sum_i E(\dot y_{i1}^2)E(\dot y_{i2}^2)E(\dot y_{i3}^2)
=n^{-2}\mu_1\mu_2\mu_3,
\end{align*}
implying
\begin{align*}
&E\left\{(r^{\otimes 3},\mu^{m-3})(r^{\otimes 2},\mu^{m-2})r_1\right\}
\\
&=n^{3/2}\left\{\mu_1\mu_4^2\cdots\mu_m^2E(\bar{\dot y}_{1110...0}\bar{\dot y}_{0110...0}\bar{\dot y}_1)+...
+\mu_1\mu_2^2\cdots\mu_{m-2}^2E(\bar{\dot y}_{10...011}\bar{\dot y}_{0...011}\bar{\dot y}_1)\right\}\\
&=n^{-1/2}(\mu_1\mu_4^2\cdots\mu_m^2\mu_1\mu_2\mu_3+...
+\mu_1\mu_2^2\cdots\mu_{m-2}^2\mu_1\mu_{m-1}\mu_m)\\
&=n^{-1/2}p_1\mu_1\sum_{1<j_1<...<j_{m-3}}\mu_{j_1}\cdots
\mu_{j_{m-3}}.
\end{align*}
Hence, via (\ref{Ed002r1}) and (\ref{Ertmkmkr1}),
\begin{equation}
E(d_{00}^2r_1)=n^{-1/2}p_1\xi_1,\label{Ed002r1b}
\end{equation}
where
\begin{align}
\xi_1&=1+\sum_{1<j_1}\mu_{j_1}+...
+\sum_{1<j_1<...<j_{m-2}}\mu_{j_1}\cdots\mu_{j_{m-2}}\notag\\
&+\mu_1\sum_{1<j_1<...<j_{m-3}}\mu_{j_1}\cdots
\mu_{j_{m-3}}.\label{xi}
\end{align}
The other terms in (\ref{Ed02r1}) may be neglected. For example, we have from (\ref{d00}) that
\begin{align*}
&E\left\{d_{00}(r^2,\mu^{m-2})r_1\right\}\\
&=E\left[\left\{\sqrt n\bar{\dot y}_{11...1}+(r^{\otimes (m-1)},\mu^1)+...+(r^{\otimes 2},\mu^{m-2})\right\}(r^2,\mu^{m-2})r_1\right],
\end{align*}
where
\begin{align*}
&E\left\{(r^{\otimes (m-k)},\mu^k)(r^2,\mu^{m-2})r_1\right\}\\
&=n^2E\left\{(\mu_1\cdots\mu_k\bar{\dot y}_{0...01...1}+...
+\mu_{m-k+1}\cdots\mu_m\bar{\dot y}_{1...10...0})\right.\\
&\left.(\mu_1\cdots\mu_{m-2}\bar{\dot y}_{m-1}\bar{\dot y}_m+...
+\mu_3\cdots\mu_m\bar{\dot y}_1\bar{\dot y}_2)\bar{\dot y}_1\right\},
\end{align*}
with e.g.
$$E(\bar{\dot y}_{1...10...0}\bar{\dot y}_1^2\bar{\dot y}_2)
=n^{-4}\sum_i\sum_j\sum_k\sum_l E(\dot y_{i1}\cdots\dot y_{i,m-k}
\dot y_{j1}\dot y_{k1}\dot y_{l2})=O(n^{-3}),$$
implying that 
$$E\left\{(r^{\otimes m-k},\mu^k)(r^{\otimes 2},\mu^{m-2})r_1\right\}
=O(n^{-1}).$$

Hence, from  (\ref{Ed002r1b}), we have
\begin{equation}
E(d_0^2r_1)=n^{-1/2}p_1\xi_1+O(n^{-1}).\label{Ed02r1b}
\end{equation}
In the following, for all $j=2,...,m$ we will define $\xi_j$ analogously as $\xi_1$.

To tackle $E(d_{0}^2C_{11})$, it follows as above that it is of the same order as $E(d_{00}^2C_{11})$. Here, from (\ref{d00}),
\begin{align}
d_{00}
&=\sqrt{n}(\bar{\dot y}_{11...1}+\mu_1\bar{\dot y}_{01...1}+...+\mu_m\bar{\dot y}_{1...10}+...\notag\\
&+\mu_1\cdots\mu_{m-2}\bar{\dot y}_{0...011}+...+\mu_3\cdots\mu_m\bar{\dot y}_{110...0}),
\label{d00b}
\end{align}
and by (\ref{C1}) and (\ref{C1b}), we may write $C_{11}=\sqrt{n}(C_1-C_{10})$ where
\begin{equation}
C_1=n^{-1}\sum_i C_{i1},\label{C1c}
\end{equation}
with
$$C_{i1}=\prod_j y_{ij}\left\{\prod_j(y_{ij}-1)-\prod_j y_{ij}\right\}.$$
Next, analogously define
\begin{equation}
X=\prod_j X_j\left\{\prod_j(X_j-1)-\prod_j X_j\right\},\label{X}
\end{equation}
where $X_j$ is a Poisson variate with parameter $\mu_j$ and all $X_j$ are simultaneously independent.

Let
\begin{align}
U_i&={\dot y}_{i1}\cdots{\dot y}_{im}+\mu_1{\dot y}_{i2}\cdots{\dot y}_{im}+...+\mu_m{\dot y}_{i1}\cdots{\dot y}_{i,m-1}+...\notag\\
&+\mu_1\cdots\mu_{m-2}{\dot y}_{i,m-1}{\dot y}_{im}+...+\mu_3\cdots\mu_m{\dot y}_{i1}{\dot y}_{i2}.\label{Ui}
\end{align}
Now, it follows from (\ref{d00b}) and (\ref{C1c}) that
$$
E(d_{00}^2C_{11})=
n^{-3/2}\sum_i\sum_j\sum_k E\{U_iU_j(C_{k1}-C_{10})\}.
$$
Here, the expectation terms are nonzero only if $i=j=k$, and so, with $\dot X=X-C_{10}$,
\begin{equation}
E(d_{00}^2C_{11})=n^{-1/2}E\left({\dot S}^2\dot X\right),\label{Ed002C11}
\end{equation}
where with $\dot X_j=X_j-\mu_j$ for all $j$,
\begin{align}
\dot S&={\dot X}_1\cdots{\dot X}_m+\mu_1{\dot X}_2\cdots{\dot X}_m+...+\mu_m{\dot X}_1\cdots{\dot X}_{m-1}+...\notag\\
&+\mu_1\cdots\mu_{m-2}{\dot X}_{m-1}{\dot X}_m+...+\mu_3\cdots\mu_m{\dot X}_1{\dot X}_2.\label{dotS}
\end{align}
The expectation in (\ref{Ed002C11}) may be worked out as a function of $\mu_1,...,\mu_m$, but we refrain to do this until later. Now, from (\ref{d11}) we get
\begin{align*}
E(d_0^2d_{11})&=p_1^{-1}\left\{C_{10}\sum_j\frac{E(d_0^2r_j)}{\mu_j}-E(d_0^2C_{11})\right\}\\
&-p_1\left\{\sum_j\mu_j^{-1}\sum_j\frac{E(d_0^2r_j)}{\mu_j}-\sum_j\frac{E(d_0^2r_j)}{\mu_j^2}\right\},
\end{align*}
where we get $E(d_0^2r_j)$ from (\ref{Ed02r1b}) and $E(d_0^2C_{11})$ from (\ref{Ed002C11}), so that
\begin{align}
E(d_0^2d_{11})&=n^{-1/2}\left\{C_{10}\sum_j\mu_j^{-1}\xi_j-p_1^{-1}E\left({\dot S}^2\dot X\right)\right.\notag\\
&\left.-p_1^2\left(s_1\sum_j\mu_j^{-1}\xi_j-\sum_j\mu_j^{-2}\xi_j\right)\right\}+O(n^{-1}),\label{Ed02d11}
\end{align}
where we have introduced
$$s_k=\sum_j\mu_j^{-k},$$
for $k=1,2,...$.

Then in turn, via (\ref{s011}) and (\ref{Ed02r1b}), we get
\begin{align}
&E(d_0^2s_{011})\notag\\
&=-\sum_j\frac{E(d_{0}^2r_j)}{\mu_j}-\frac{E(d_0^2d_{11})}{d_{10}}\notag\\
&=-n^{-1/2}\left[p_1\sum_j\mu_j^{-1}\xi_j
+d_{10}^{-1}\left\{C_{10}\sum_j\mu_j^{-1}\xi_j-p_1^{-1}E\left({\dot S}^2\dot X\right)\right.\right.\notag\\
&\left.\left.-p_1^2\left(s_1\sum_j\mu_j^{-1}\xi_j-\sum_j\mu_j^{-2}\xi_j\right)\right\}\right]+O(n^{-1}).
\label{Ed02s011}
\end{align}

To go on with $E(d_0^2s_{012})$, we at first note that, from (\ref{d00b}),
$$
E(d_{00}^2r_1^2)=
n^{-2}\sum_i\sum_j\sum_k\sum_l E\{U_iU_j{\dot y}_{k1}{\dot y}_{l1}\}.
$$
The expectation inside of the sum is nonzero not only for $i=j=k=l$ (which creates a $O(n^{-1})$ term) but also for $i=j$ and $k=l$, where $i$ and $k$ are not necessarily equal. It follows as above that
\begin{equation}
E(d_{00}^2r_1^2)=E(\dot S^2)E(\dot X_1^2)+O(n^{-1})=\mu_1E(\dot S^2)+O(n^{-1}).\label{Ed002r1r1}
\end{equation}
Similarly, it is seen that $E(d_{00}^2r_1r_2)=O(n^{-1})$, since in this case, such a 'cross product term' cannot arise.

Next, to calculate $E(d_{00}^2d_{11}r_1)$, we get as above that, with $\dot X=X-C_{10}$ (observe that $E(\dot S X_j)=0$ for all $j$),
\begin{align}
&E(d_{00}^2C_{11}r_1)\notag\\&=
n^{-2}\sum_i\sum_j\sum_k\sum_l E\{U_iU_j(C_{k1}-C_{10}){\dot y}_{l1}\}
\notag\\
&=E(\dot S^2)E(\dot X\dot X_1)
+O(n^{-1}).\label{Ed002C11r1}
\end{align}
Thus, via (\ref{d11}) and (\ref{Ed002r1r1}),
\begin{align}
&E(d_{00}^2d_{11}r_1)\notag\\
&=p_1^{-1}\left\{C_{10}\frac{E(d_{00}^2r_1^2)}{\mu_1}-E(d_{00}^2C_{11}r_1)\right\}\notag\\
&-p_1\left\{\sum_j\mu_j^{-1}\frac{E(d_{00}^2r_1^2)}{\mu_1}
-\frac{E(d_{00}^2r_1^2)}{\mu_1^2}\right\}\notag\\
&=p_1^{-1}E({\dot S}^2)\left\{C_{10}-E(\dot X\dot X_1)\right\}-p_1E({\dot S}^2)(s_1-\mu_1^{-1})
\notag\\&+O(n^{-1}).\label{Ed002d11r1}
\end{align}

Next, (\ref{d11}) and simplifications yield
\begin{align*}
d_{11}^2&=
p_1^{-2}C_{11}^2+2DC_{11}\sum_j\frac{r_j}{\mu_j}+p_1^{2}D^2\left(\sum_j\frac{r_j}{\mu_j}\right)^2\\
&-2C_{11}\sum_j\frac{r_j}{\mu_j^2}+p_1^2\left(\sum_j\frac{r_j}{\mu_j^2}\right)^2-2p_1^2D\sum_j\frac{r_j}{\mu_j}\sum_j\frac{r_j}{\mu_j^2},
\end{align*}
where
\begin{equation}
D=s_1-p_1^{-2}C_{10}
=p_1^{-1}d_{10}+2s_1,\label{D}
\end{equation}
following from the fact that, via (\ref{md1}) and (\ref{d10}),
\begin{equation}
-d_{10}=p_1^{-1}C_{10}+p_1s_1.
\label{md10}
\end{equation} 
To find $E(d_{00}^2d_{11}^2)$, we need
\begin{align}
E(d_{00}^2C_{11}^2)&=
n^{-2}\sum_i\sum_j\sum_k\sum_l E\{U_iU_j(C_{k1}-C_{10})(C_{l1}-C_{10})\}
\notag\\
&=E(\dot S^2)E({\dot X}^2)+2\left\{E(\dot S\dot X)\right\}^2
+O(n^{-1}).\label{Ed00C112}
\end{align}
Now, together with (\ref{Ed002C11r1}) and (\ref{Ed002r1r1}), we get
\begin{align}
&E(d_{00}^2d_{11}^2)\notag\\
&=p_1^{-2}E(d_{00}^2C_{11}^2)+2D\sum_j\mu_j^{-1}E(d_{00}^2C_{11}r_j)+p_1^2D^2\sum_j\mu_j^{-2}E(d_{00}^2r_j^2)\notag\\
&-2\sum_j\mu_j^{-2}E(d_{00}^2C_{11}r_j)+p_1^2\sum_j\mu_j^{-4}E(d_{00}^2r_j^2)-2p_1^2D\sum_j\mu_j^{-3}E(d_{00}^2r_j^2)\notag\\
&=2p_1^{-2}\left\{E(\dot S\dot X)\right\}^2+E(\dot S^2)
\left\{p_1^{-2}E({\dot X}^2)
+2DA_1+p_1^2D^2s_1\right.\notag\\
&\left.
-2A_2+p_1^2s_3-2p_1^2Ds_2\right\}+O(n^{-1/2}),\label{Ed002d112}
\end{align}
where
$$A_1=\sum_j\mu_j^{-1}E(\dot X\dot X_j),\quad A_2=\sum_j\mu_j^{-2}E(\dot X\dot X_j).$$
Similarly, from (\ref{d12}),  (\ref{Ed002r1r1}) and (\ref{Ed002C11r1}),
\begin{align}
&E(d_{00}^2d_{12})\notag\\&=
p_1^{-1}\left\{\sum_j\mu_j^{-1}E(d_{00}^2C_{11}r_j)-C_{10}
\sum_j\mu_j^{-2}E(d_{00}^2r_j^2)\right\}+O(n^{-1})\notag\\
&=p_1^{-1}E(\dot S^2)\left(A_1-C_{10}s_1\right)+O(n^{-1}).\label{Ed002d12}
\end{align}
Next, (\ref{s012}) implies
\begin{align*}
E(d_{00}^2s_{012})&=\sum_j\mu_j^{-2}E(d_{00}^2r_j^2)
+d_{10}^{-1}\sum_j\mu_j^{-1}E(d_{00}^2d_{11}r_j)\\
&+d_{10}^{-2}E(d_{00}^2d_{11}^2)-d_{10}^{-1}E(d_{00}^2d_{12})
+O(n^{-1}),
\end{align*}
and via (\ref{Ed002r1r1}), (\ref{Ed002d11r1}), (\ref{D}), (\ref{Ed002d112}), (\ref{Ed002d12}) and simplifications,
\begin{equation}
E(d_{00}^2s_{012})=2d_{10}^{-2}p_1^{-2}\left\{E(\dot S\dot X)\right\}^2+E({\dot S}^2)F+O(n^{-1}),\label{Ed002s012}\end{equation}
with
\begin{align}
F&=p_1^{-2}d_{10}^{-2}E({\dot X}^2)+4 d_{10}^{-2}s_1A_1
-2d_{10}^{-2}A_2+p_1d_{10}^{-1}(s_1^2-s_2)\notag\\
&+p_1^2d_{10}^{-2}(4s_1^3-4s_1s_2+s_3).\label{F}
\end{align}
Now, from (\ref{f01b}), we get
\begin{equation}
E(f_{01})=p_1^{-1}d_{10}^{-1}\left\{E(d_{00}^2)
+n^{-1/2}E(d_{00}^2s_{011})+n^{-1}E(d_{00}^2s_{012})\right\},
\label{Ef01}
\end{equation}
where we get $E(d_{00}^2)$ from (\ref{Ed02b}), $E(d_{00}^2s_{011})$ from (\ref{Ed02s011}) and $E(d_{00}^2s_{012})$ from (\ref{Ed002s012}).

Next, (\ref{f02}), (\ref{Pinv}) and (\ref{md1inv}) imply
\begin{equation}
f_{02}=\frac{d_0^2d_{32}}{2p_1^2d_{10}^2}G,\label{f02b}
\end{equation}
where
\begin{align}
G=&\left(1-n^{-1/2}\sum_j\frac{r_j}{\mu_j}+n^{-1}\sum_{j\geq k}\frac{r_jr_k}{\mu_j\mu_k}\right)^2\notag\\
&\cdot\left(1-n^{-1/2}\frac{d_{11}}{d_{10}}+n^{-1}\frac{d_{11}^2-d_{10}d_{12}}{d_{10}^2}\right)^2+O(n^{-3/2})\notag\\
=&1-2n^{-1/2}G_1+n^{-1}G_2+O(n^{-3/2}),\label{G}
\end{align}
with
\begin{align}
G_1&=\sum_j\frac{r_j}{\mu_j}+\frac{d_{11}}{d_{10}},\label{G1}\\
G_2&=2\sum_{j\geq k}\frac{r_jr_k}{\mu_j\mu_k}+\left(\sum_j\frac{r_j}{\mu_j}\right)^2+4\frac{d_{11}}{d_{10}}\sum_j\frac{r_j}{\mu_j}+2\frac{d_{11}^2-d_{10}d_{12}}{d_{10}^2}+\frac{d_{11}^2}{d_{10}^2}.
\label{G2}
\end{align}
We get $d_{32}$ from (\ref{d32}), where (cf (\ref{d0b}))
\begin{align}
&2\bar y_{11...1}-P=2\bar y_{11...1}-\prod_j\bar y_j\notag\\
&=2n^{-1}\sum_i(\mu_1+\dot y_{i1})\cdots(\mu_m+\dot y_{im})-(\mu_1+\bar{\dot y}_1)\cdots(\mu_m+\bar{\dot y}_m)\notag\\
&=2p_1+2(\bar{\dot y}_{11...1}+\mu_1\bar{\dot y}_{01...1}+...+\mu_m\bar{\dot y}_{1...10}
+...\notag\\
&+\mu_1\mu_2\cdots\mu_{m-2}\bar{\dot y}_{0...011}+...+\mu_3\cdots\mu_m\bar{\dot y}_{110...0})\notag\\
&+2(\mu_1\cdots\mu_{m-1}\bar{\dot y}_m+...+\mu_2\cdots\mu_m\bar{\dot y}_1)\notag\\
&-p_1-(\mu_1\cdots\mu_{m-1}\bar{\dot y}_m+...+\mu_2\cdots\mu_m\bar{\dot y}_1)\notag\\
&-(\mu_1\cdots\mu_{m-2}\bar{\dot y}_{m-1}\bar{\dot y}_m+...+\mu_3\cdots\mu_m\bar{\dot y}_1\bar{\dot y}_2)-...\notag\\
&=p_1+n^{-1/2}\left\{2d_{00}+(r,\mu^{m-1})\right\}
-n^{-1}(r^2,\mu^{m-2})+O(n^{-3/2}),\label{d32term1}
\end{align} 
and
\begin{align}
&\by(1^{m-1},0^1)\notag\\
&=\bar y_{01...1}+...+\bar y_{1...10}\notag\\
&=n^{-1}\sum_i(\mu_2+\dot y_{i2})\cdots(\mu_m+\dot y_{im})+...\notag\\
&+n^{-1}\sum_i(\mu_1+\dot y_{i1})\cdots(\mu_{m-1}+\dot y_{i,m-1})
\notag\\
&=\mu_2\cdots\mu_m+\mu_2\cdots\mu_{m-1}\bar{\dot y}_m+...+\mu_3\cdots\mu_m\bar{\dot y}_2+...\notag\\
&+\mu_2\bar{\dot y}_{00...11}+...+\mu_m\bar{\dot y}_{01...10}+\bar{\dot y}_{01...1}+...\notag\\
&+\mu_1\cdots\mu_{m-1}+\mu_1\cdots\mu_{m-2}\bar{\dot y}_{m-1}+...+\mu_2\cdots\mu_{m-1}\bar{\dot y}_1+...\notag\\
&+\mu_1\bar{\dot y}_{01...10}+...+\mu_{m-1}\bar{\dot y}_{1...100}+\bar{\dot y}_{1...10}\notag\\
&=p_1\sum_j\mu_j^{-1}+n^{-1/2}H,\label{d32term2}
\end{align}
where, with $\dot r_{12}=\sqrt{n}\bar{\dot y}_{110...0}$ ectetera,
\begin{align}
H&=\sum_{j_1<...<j_{m-2}}
\mu_{j_1}\cdots\mu_{j_{m-2}}\sum_{k\notin(j_1,...,j_{m-2})}r_k+...
\notag\\
&+\sum_j\mu_j\sum_{k_1<...<k_{m-2},\forall k_i\neq j}\dot r_{k_1,...,k_{m-2}}
+\sum_j\sum_{k_1<...<k_{m-1},\forall k_i\neq j}\dot r_{k_1,...,k_{m-1}}.\label{H}
\end{align}
Moreover, since by (\ref{e2bar}) and (\ref{C1b}),
$$\bar e_2=C_1=C_{10}+n^{-1/2}C_{11},$$
we get from (\ref{d32}), (\ref{d32term1}) and (\ref{d32term2}) that
\begin{align}
d_{32}&=
\left[p_1+n^{-1/2}\left\{2d_{00}+(r,\mu^{m-1})\right\}
-n^{-1}(r^2,\mu^{m-2})\right]\notag\\
&\cdot\left(p_1s_1+n^{-1/2}H\right)
+C_{10}+n^{-1/2}C_{11}+O(n^{-3/2})\notag\\
&=d_{320}+n^{-1/2}d_{321}+n^{-1}d_{322}+O(n^{-3/2}),\label{d32b}
\end{align}
where
\begin{align}
d_{320}&=p_1^2s_1+C_{10}
,\label{d320}\\
d_{321}&=p_1H+p_1s_1
\left\{2d_{00}+(r,\mu^{m-1})\right\}+C_{11}
,\label{d321}\\
d_{322}&=H\left\{2d_{00}+(r,\mu^{m-1})\right\}-
p_1s_1(r^2,\mu^{m-2})
.\label{d322}
\end{align}
Hence, inserting (\ref{G}) and (\ref{d32b}) into (\ref{f02b}), we get
\begin{equation}
f_{02}=\frac{d_0^2}{2p_1^2d_{10}^2}(d_{320}+n^{-1/2}s_{021}+n^{-1}s_{022})+O(n^{-3/2}),\label{f02c}
\end{equation}
where
\begin{align}
s_{021}&=d_{321}-2d_{320}G_1,\label{s021}\\
s_{022}&=d_{322}-2d_{321}G_1+d_{320}G_2.\label{s022}
\end{align}
Now, observing that from (\ref{d320}) and (\ref{md10}), 
\begin{equation}
d_{320}=-p_1d_{10},\label{d320b}
\end{equation} 
we get via (\ref{Ed02b}) that
\begin{align}
E(d_0^2d_{320})&=-p_1^2d_{10}^2\left(1-n^{-1}d_{10}^{-1}\sum_{j_1<...<j_{m-2}}\mu_{j_1}\cdots\mu_{j_{m-2}}\right)\notag\\&+O(n^{-3/2}).\label{Ed02d320}
\end{align}
To tackle $E(d_0^2d_{321})$, we at first need $E(d_0^2H)$. To this end, we have from  (\ref{d00b}), (\ref{Ui}) and (\ref{dotS}) that
$$
E(d_{00}^2H)=n^{-3/2}\sum_i\sum_j\sum_lE(U_iU_jH_l)=n^{-1/2}E(\dot S^2Z),
$$
where, from (\ref{H}),
\begin{align}
H_l&=\sum_{j_1<...<j_{m-2}}
\mu_{j_1}\cdots\mu_{j_{m-2}}\sum_{k\notin(j_1,...,j_{m-2})}\dot y_{lk}+...\notag\\
&+\sum_j\mu_j\sum_{k_1<...<k_{m-2},\forall k_i\neq j}\dot y_{l,k_1}\cdots\dot y_{l,k_{m-2}}\notag\\
&+\sum_j\sum_{k_1<...<k_{m-1},\forall k_i\neq j}\dot y_{l,k_1}\cdots\dot y_{l,k_{m-1}},\label{Hl}
\end{align}
and
\begin{align}
Z\notag&=\sum_{j_1<...<j_{m-2}}
\mu_{j_1}\cdots\mu_{j_{m-2}}\sum_{k\notin(j_1,...,j_{m-2})}\dot X_k+...
\notag\\
&+\sum_j\mu_j\sum_{k_1<...<k_{m-2},\forall k_i\neq j}\dot X_{k_1}\cdots\dot X_{k_{m-2}}\notag\\&+\sum_j\sum_{k_1<...<k_{m-1},\forall k_i\neq j}\dot X_{k_1}\cdots\dot X_{k_{m-1}}.\label{Z}
\end{align}
Via (\ref{d0b}), it follows that $E(d_0^2H)$ asymptotically behaves like $E(d_{00}^2H)$, i.e.
\begin{equation}
E(d_0^2H)=n^{-1/2}E(\dot S^2Z)+O(n^{-1}).\label{Ed02H}
\end{equation}

Moreover, from (\ref{d0b}),
\begin{equation}
E(d_0^2d_{00})=E(d_{00}^3)-2n^{-1/2}E\left\{d_{00}^2(r^2,\mu^{m-2})\right\}+O(n^{-1}),\label{Ed02d00}
\end{equation}
where via (\ref{d00}),
\begin{align}
&E(d_{00}^3)\notag\\
&=n^{3/2}E(\bar{\dot y}_{11...1}^3)+E\left\{(r^{\otimes(m-1)},\mu^1)^3\right\}+...+E\left\{(r^{\otimes 2},\mu^{m-2})^3\right\}
\notag\\
&+3nE\left\{\bar{\dot y}_{11...1}^2(r^{\otimes(m-1)},\mu^1)\right\}
+...+3E\left\{(r^{\otimes 3},\mu^{m-3})(r^{\otimes 2},\mu^{m-2})^2\right\}.
\label{Ed003}
\end{align}
Here, as before, the 'cross product terms' are zero, while
\begin{align}
&n^{3/2}E(\bar{\dot y}_{11...1}^3)\notag\\
&=n^{-3/2}E\left\{\left(\sum_i \dot y_{i1}\cdots \dot y_{im}\right)^3\right\}\notag\\
&=n^{-3/2}\sum_i\sum_j\sum_k E(\dot y_{i1}\cdots \dot y_{im}\dot y_{j1}\cdots \dot y_{jm}\dot y_{k1}\cdots \dot y_{km})\notag\\
&=n^{-3/2}\sum_i E(\dot y_{i1}^3)\cdots E(\dot y_{im}^3)=
n^{-1/2}\mu_1\cdots\mu_m=n^{-1/2}p_1\label{Ebary3}
\end{align}
and
\begin{align*}
&E\left\{(r^{\otimes(m-k)},\mu^k)^3\right\}\\
&=n^{3/2}E\left\{(\mu_1\mu_2\cdots\mu_k\bar{\dot y}_{0...01...1}+...+\mu_{m-k+1}\cdots\mu_m\bar{\dot y}_{1...10...0})^3\right\}\\
&=n^{3/2}\left\{\mu_1^3\cdots \mu_k^3E(\bar{\dot y}_{0...01...1}^3)
+...+\mu_{m-k+1}^3\cdots \mu_m^3E(\bar{\dot y}_{1...10...0}^3)\right\},
\end{align*}
where as in (\ref{Ebary3}),
$$n^{3/2}E(\bar{\dot y}_{1...10...0}^3)=n^{-1/2}\mu_1\cdots\mu_{m-k}$$
etcetera, and it follows that
\begin{equation}
E\left\{(r^{\otimes(m-k)},\mu^k)^3\right\}=n^{-1/2}p_1\sum_{j_1<...<j_k}\mu_{j_1}\cdots\mu_{j_k}.\label{Ermu3}
\end{equation}
Conclusively, inserting (\ref{Ebary3}) and (\ref{Ermu3}) into (\ref{Ed003}), we find from (\ref{d10}) that
\begin{equation}
E(d_{00}^3)=n^{-1/2}p_1d_{10}.\label{Ed003b}
\end{equation}

Next, (\ref{d00}) implies
\begin{align}
&E\left\{d_{00}^2(r^2,\mu^{m-2})\right\}\notag\\&=
E\left[\left\{\sqrt{n}\bar{\dot y}_{11...1}+(r^{\otimes m-1},\mu^1)+...+(r^{\otimes 2},\mu^{m-2})\right\}^2(r^2,\mu^{m-2})\right]\notag\\
&=nE\left\{\bar{\dot y}_{11...1}^2(r^2,\mu^{m-2})\right\}
+E\left\{(r^{\otimes m-1},\mu^1)^2(r^2,\mu^{m-2})\right\}+...\notag\\
&+E\left\{(r^{\otimes 2},\mu^{m-2})^2(r^2,\mu^{m-2})\right\},
\label{Ed002r2mum2}
\end{align}
since 'cross terms' have expectation zero. Here,
\begin{align*}
E\left\{\bar{\dot y}_{11...1}^2(r^2,\mu^{m-2})\right\}
&=E\left\{\bar{\dot y}_{11...1}^2(\mu_1\cdots\mu_{m-2}r_{m-1}r_m+...
+\mu_3\cdots\mu_mr_1r_2)\right\}\\
&=\mu_1\cdots\mu_{m-2}E(\bar{\dot y}_{11...1}^2r_{m-1}r_m)+...
+\mu_3\cdots\mu_mE(\bar{\dot y}_{11...1}^2r_1r_2).
\end{align*}
Here,
$$nE(\bar{\dot y}_{11...1}^2r_1r_2)=n^{-2}\sum_i\sum_j\sum_k\sum_l E(\dot y_{i1}\cdots\dot y_{im}\dot y_{j1}\cdots\dot y_{jm}\dot y_{k1}\dot y_{l2})=O(n^{-1})$$
(in fact, this is zero for $m>2$), and it is similarly seen that all other terms of (\ref{Ed002r2mum2}) are also $O(n^{-1})$. Hence, $E\left\{d_{00}^2(r^2,\mu^{m-2})\right\}=O(n^{-1})$, and so, via (\ref{Ed02d00}) and (\ref{Ed003b}), we have
\begin{equation}
E(d_0^2d_{00})=n^{-1/2}p_1d_{10}+O(n^{-3/2}).\label{Ed02d00b}
\end{equation}
Furthermore, we similarly have
\begin{align}
&E\left\{d_{00}^2(r,\mu^{m-1})\right\}\notag\\&=
E\left[\left\{\sqrt{n}\bar{\dot y}_{11...1}+(r^{\otimes m-1},\mu^1)+...+(r^{\otimes 2},\mu^{m-2})\right\}^2(r,\mu^{m-1})\right]\notag\\
&=nE\left\{\bar{\dot y}_{11...1}^2(r,\mu^{m-1})\right\}
+E\left\{(r^{\otimes m-1},\mu^1)^2(r,\mu^{m-1})\right\}+...\notag\\
&+E\left\{(r^{\otimes 2},\mu^{m-2})^2(r,\mu^{m-1})\right\},
\label{Ed002rmum1}
\end{align}
where
\begin{align*}
E\left\{\bar{\dot y}_{11...1}^2(r,\mu^{m-1})\right\}
&=E\left\{\bar{\dot y}_{11...1}^2(\mu_1\cdots\mu_{m-1}r_m+...
+\mu_2\cdots\mu_mr_1)\right\}\\
&=\mu_1\cdots\mu_{m-1}E(\bar{\dot y}_{11...1}^2r_m)+...
+\mu_2\cdots\mu_mE(\bar{\dot y}_{11...1}^2r_1),
\end{align*}
where (\ref{Ebary2r1}) yields
$$nE(\bar{\dot y}_{11...1}^2r_1)=n^{-1/2}p_1.$$
It follows that
\begin{equation}
nE\left\{\bar{\dot y}_{11...1}^2(r,\mu^{m-1})\right\}
=n^{-1/2}p_1\sum_{j_1<...<j_{m-1}}\mu_{j_1}\cdots\mu_{j_{m-1}}.\label{Ey112rmum1}
\end{equation}
Moreover, for $k=1,...,m-2$, we have
\begin{align*}
&E\left\{(r^{\otimes m-k},\mu^k)^2(r,\mu^{m-1})\right\}\\
&=n^{3/2}E\left\{(\mu_1\cdots\mu_k\bar{\dot y}_{0...01...1}+...
+\mu_{m-k+1}\cdots\mu_m\bar{\dot y}_{1...10...0})^2\right.\\
&\left.(\mu_1\cdots\mu_{m-1}\bar{\dot y}_m+...+\mu_2\cdots\mu_m\bar{\dot y}_1)\right\}
\end{align*}
where from (\ref{Ebary2bary1}),
$$n^{3/2}E(\bar{\dot y}_{1...10...0}^2\bar{\dot y}_1)
=n^{-1/2}\mu_1\cdots\mu_{m-k}.$$
This type of results imply
\begin{align*}
&E\left\{(r^{\otimes m-k},\mu^k)^2(r,\mu^{m-1})\right\}\\
&=n^{-1/2}\left(\mu_2\cdots\mu_m\sum_{1<l_1<...<l_k}\mu_{l_1}^2\cdots\mu_{l_k}^2\prod_{k\notin(l_1,...,l_k)}\mu_k\right.
\\&\left.+...+\mu_1\cdots\mu_{m-1}\sum_{l_1<...<l_k<m}\mu_{l_1}^2\cdots\mu_{l_k}^2\prod_{k\notin(l_1,...,l_k)}\mu_k\right)\\
&=n^{-1/2}p_1\sum_{j_1<...<j_{m-1}}\mu_{j_1}\cdots\mu_{j_{m-1}}
\sum_{l_1<...<l_k,\forall l_j\in(j_1,...,j_{m-1})}\mu_{l_1}\cdots\mu_{l_k}.
\end{align*}
Hence, via (\ref{Ed002rmum1}) and (\ref{Ey112rmum1}),
\begin{align*}
&E\left\{d_{00}^2(r,\mu^{m-1})\right\}\\&=
n^{-1/2}p_1\sum_{j_1<...<j_{m-1}}\mu_{j_1}\cdots\mu_{j_{m-1}}
\left(1+\sum_{k=1}^{m-2}\sum_{l_1<...<l_k,\forall l_j\in(j_1,...,j_{m-1})}\mu_{l_1}\cdots\mu_{l_k}\right)\\
&=n^{-1/2}\eta.
\end{align*}
Here, an alternative expression more suitable for future calculations is
\begin{equation}
\eta=p_1^2\sum_j\mu_j^{-1}\left\{\prod_{k\neq j}(1+\mu_k)-\prod_{k\neq j}\mu_k\right\}.\label{eta}
\end{equation}

In view of (\ref{d0b}), this also means that 
\begin{equation}
E\left\{d_0^2(r,\mu^{m-1})\right\}=n^{-1/2}\eta+O(n^{-1}),\label{Ed002rmum1b}
\end{equation}
because in the usual manner, it can be seen that
$$E\{d_{00}(r^2,\mu^{m-2})(r,\mu^{m-1})\}=O(n^{-1}).$$
Now, from (\ref{d321}),
\begin{align*}
&E(d_0^2d_{321})\\
&=p_1E(d_0^2H)+p_1s_1\left[2E(d_0^2d_{00})+E\{d_0^2(r,\mu^{m-1})\}\right]+E(d_0^2C_{11}),
\end{align*}
and via (\ref{Ed02H}), (\ref{Ed02d00b}), (\ref{Ed002rmum1b}), (\ref{Ed002C11}) and (\ref{d0b}), we get
\begin{align}
E(d_0^2d_{321})&=n^{-1/2}\left\{p_1E(\dot S^2Z)\right.\notag\\
&\left.+p_1s_1(2p_1d_{10}+\eta
)+E(\dot S^2\dot X)\right\}+O(n^{-1}).\label{Ed02d321}
\end{align}
Next, via (\ref{G1}), (\ref{Ed02r1b}) and (\ref{Ed02d11}), we have
\begin{align*}
&E(d_0^2G_1)\\
&=\sum_j\mu_j^{-1}E(d_0^2r_j)+d_{10}^{-1}E(d_0^2d_{11})\\
&=n^{-1/2}\left[p_1\sum_j\mu_j^{-1}\xi_j
+d_{10}^{-1}
\left\{C_{10}\sum_j\mu_j^{-1}\xi_j-p_1^{-1}E\left({\dot S}^2\dot X\right)\right.\right.\\
&\left.\left.-p_1^2\left(s_1\sum_j\mu_j^{-1}\xi_j-\sum_j\mu_j^{-2}\xi_j\right)\right\}\right]+O(n^{-1}).
\end{align*}
Here, from (\ref{md10}), we have 
$$C_{10}=-p_1d_{10}-p_1^2s_1,$$
which yields
\begin{align*}
&E(d_0^2G_1)\\
&=-n^{-1/2}d_{10}^{-1}\left\{p_1^2\left(2s_1\sum_j\mu_j^{-1}\xi_j-\sum_j\mu_j^{-2}\xi_j\right)
+p_1^{-1}E\left(\dot S^2\dot X\right)\right\}\\
&+O(n^{-1}).
\end{align*}
Hence, together with (\ref{Ed02d321}), (\ref{d321}) and (\ref{d320b}), we have
\begin{align}
&E(d_0^2s_{021})\notag\\
&=E(d_0^2d_{321})+2p_1d_{10}E(d_0^2G_1)\notag\\
&=n^{-1/2}\left\{p_1E(\dot S^2Z)+p_1s_1(2p_1d_{10}+\eta)
+E(\dot S^2\dot X)\right.\notag\\
&\left.-2p_1^3\left(2s_1\sum_j\mu_j^{-1}\xi_j-\sum_j\mu_j^{-2}\xi_j\right)-2E(\dot S^2\dot X)
\right\}+O(n^{-1}).\label{Ed02s021}
\end{align}

To find $E(d_0^2s_{022})$, then in view of (\ref{s022}), we at first need $E(d_0^2d_{322})$. In view of (\ref{d322}), we start with $E(d_0^3H)$.  From (\ref{d00b}), (\ref{Ui}), (\ref{dotS}), (\ref{Hl}) and (\ref{Z}), 
$$E(d_{00}^3H)=n^{-2}\sum_i\sum_j\sum_k\sum_l E(U_iU_jU_kH_l)=3E(\dot S^2)E(\dot S Z)+O(n^{-1}),$$
and because by (\ref{d0b}), we have
\begin{equation}
E(d_0^2d_{00}H)=3E(\dot S^2)E(\dot S Z)+O(n^{-1}).\label{Ed02d00H}
\end{equation}
Similarly, we get
\begin{align*}
&E\{d_{00}^2H(r,\mu^{m-1})\}\\
&=n^{-2}\sum_i\sum_j\sum_k\sum_l E\{U_iU_jH_k(\mu_1\cdots\mu_{m-1}\dot y_{lm}+...+\mu_2\cdots\mu_m\dot y_{l1})\}\\
&=E(\dot S^2)E(ZV),
\end{align*}
where
\begin{equation}
V=\mu_1\cdots\mu_{m-1}\dot X_m+...+\mu_2\cdots\mu_m\dot X_1,
\label{V}
\end{equation}
implying
\begin{equation}
E\{d_0^2H(r,\mu^{m-1})\}=E(\dot S^2)E(ZV)+O(n^{-1}).
\label{Ed02Hrmu}
\end{equation}
Moreover,
\begin{align*}
&E\{d_{00}^2(r^2,\mu^{m-2})\}\\
&=n^{-2}\sum_i\sum_j\sum_k\sum_l E\{U_iU_j(\mu_1\cdots\mu_{m-2}\dot y_{k,m-1}\dot y_{l,m}+...+\mu_2\cdots\mu_m\dot y_{k1}\dot y_{l2})\}\\
&=O(n^{-1}),
\end{align*}
implying $E\{d_0^2(r^2,\mu^{m-2})\}=O(n^{-1})$. Hence, via (\ref{d322}), (\ref{Ed02d00H}) and (\ref{Ed02Hrmu}), we get
\begin{equation}
E(d_0^2d_{322})=E(\dot S^2)\{6E(\dot S Z)+E(ZV)\}+O(n^{-1}).
\label{Ed02d332}
\end{equation}

As for $E(d_0^2d_{321}G_1)$, we start with $E(d_0^2d_{321}r_1)$, and to this end, we have in the usual manner
\begin{align} 
E(d_{00}^2Hr_1)&=n^{-2}\sum_i\sum_j\sum_k\sum_lE(U_iU_jH_k\dot y_{l1})\notag\\
&=E(\dot S^2)E(Z\dot X_1)+O(n^{-1}),\label{Ed002Hr1}
\end{align}
\begin{equation}
E(d_{00}^3r_1)=n^{-2}\sum_i\sum_j\sum_k\sum_lE(U_iU_jU_k\dot y_{l1})=0,\label{Ed003r1}
\end{equation}
and
\begin{align}
&E\{d_{00}^2(r,\mu^{m-1})r_1\}\notag\\
&=n^{-2}\sum_i\sum_j\sum_k\sum_l
E\{U_iU_j(\mu_1\cdots\mu_{m-1}\dot y_{km}+...+\mu_2\cdots\mu_m\dot y_{m1})\dot y_{l1}
\}\notag\\
&=p_1E(\dot S^2)+O(n^{-1}),\label{Ed002rmum1r1}
\end{align}
which together with (\ref{Ed002C11r1}) and (\ref{d321}) yield
\begin{align}
&E(d_0^2d_{321}r_1)\notag\\
&=E(\dot S^2)\{p_1E(Z\dot X_1)+p_1^2s_1+E(\dot X\dot X_1)\}
+O(n^{-1})\notag\\
&=E(\dot S^2)\{p_1E(\tilde Z\dot X_1)+p_1^2s_1\},\label{Ed02d321r1}
\end{align}
where $\tilde Z=Z+p_1^{-1}\dot X$.
Similarly, to find $E(d_0^2d_{321}d_{11})$, we need, by (\ref{d11}),
\begin{align}
&E(d_{00}^2Hd_{11})\notag\\&=
p_1^{-1}\left\{C_{10}\sum_j\mu_j^{-1}E(d_{00}^2Hr_j)-E(d_{00}^2HC_{11})\right\}\notag\\
&-p_1\left\{s_1\sum_j\mu_j^{-1}E(d_{00}^2Hr_j)
-\sum_j\mu_j^{-2}E(d_{00}^2Hr_j)\right\}.\label{Ed002Hd11}
\end{align}
Here, we have in the usual manner that
\begin{align} 
E(d_{00}^2HC_{11})&=n^{-2}\sum_i\sum_j\sum_k\sum_lE\{U_iU_jH_k\dot (C_{l1}-C_{10})\}\notag\\
&=E(\dot S^2)E(Z\dot X)+2E(\dot S Z)E(\dot S\dot X)+O(n^{-1}),\label{Ed002HC11}
\end{align}
implying via (\ref{Ed002Hr1}) and (\ref{Ed002Hd11}) that
\begin{align}
&E(d_0^2Hd_{11})\notag\\&=
p_1^{-1}E(\dot S^2)\left\{C_{10}\sum_j\mu_j^{-1}E(Z\dot X_j)-E(Z\dot X)\right\}\notag\\
&-2p_1^{-1}E(\dot S Z)E(\dot S\dot X)\notag\\
&-p_1E(\dot S^2)\left\{s_1\sum_j\mu_j^{-1}E(Z\dot X_j)-\sum_j\mu_j^{-2}E(Z\dot X_j)\right\}\notag\\
&+O(n^{-1}).\label{Ed02Hd11}
\end{align}
Similarly, via (\ref{d11}) and (\ref{Ed003r1}), we have
$$E(d_{00}^3d_{11})=-p_1^{-1}E(d_{00}^3C_{11}),$$
where
\begin{align}
E(d_{00}^3C_{11})&=n^{-2}\sum_i\sum_j\sum_k\sum_lE\{U_iU_jU_k (C_{l1}-C_{10})\}\notag\\
&=3E(\dot S^2)E(\dot S\dot X)+O(n^{-1}),\label{Ed003C11}
\end{align}
and so,
\begin{equation}
E(d_0^2d_{00}d_{11})=-3p_1^{-1}E(\dot S^2)E(\dot S\dot X)+O(n^{-1}).\label{Ed03d11}
\end{equation}
To find $E\{d_{00}^2(r,\mu^{m-1})d_{11}\}$, we need (cf (\ref{V}) and (\ref{Ed02Hrmu}))
\begin{align*}
&E\{d_{00}^2(r,\mu^{m-1})C_{11}\}\\
&=n^{-2}\sum_i\sum_j\sum_k\sum_l E\{U_iU_j(\mu_1\cdots\mu_{m-1}\dot y_{km}+...+\mu_2\cdots\mu_m\dot y_{k1})(C_{l1}-C_{10})\}\\
&=E(\dot S^2)E(V\dot X)+O(n^{-1}),
\end{align*}
Together with (\ref{d11}) and (\ref{Ed002rmum1r1}), this gives
\begin{align*}
E\{d_{00}^2(r,\mu^{m-1})d_{11}\}
&=p_1^{-1}E(\dot S^2)
\left\{C_{10}p_1s_1-E(V\dot X)\right\}\\
&-p_1^2E(\dot S^2)(s_1^2-s_2)+O(n^{-1}).
\end{align*}
and because from (\ref{md10}), $C_{10}=-p_1d_{10}-p_1^2s_1$ which implies
\begin{align}
&E\{d_{00}^2(r,\mu^{m-1})d_{11}\}\notag\\
&=E(\dot S^2)\left\{-p_1d_{10}s_1-p_1^{-1}E(V\dot X)-2p_1^2s_1^2+p_1^2s_2\right\}\notag\\
&+O(n^{-1}).\label{Ed002rmum1d11}
\end{align}
As for $E(d_{00}^2C_{11}d_{11})$, we have
\begin{align}
E(d_{00}^2C_{11}^2)&=n^{-2}\sum_i\sum_j\sum_k\sum_lE\{U_iU_j(C_{k1}-C_{10})(C_{l1}-C_{10})\}\notag\\
&=E(\dot S^2)E(\dot X^2)+2\{E(\dot S\dot X)\}^2+O(n^{-1}),
\label{Ed002C112}
\end{align}
implying via (\ref{Ed002C11r1}) and (\ref{d11}) that
\begin{align*}
&E(d_{00}^2C_{11}d_{11})\\
&=p_1^{-1}E(\dot S^2)\left\{C_{10}\sum_j\mu_j^{-1}E(\dot X\dot X_j)-E(\dot X^2)\right\}-2p_1^{-1}\left\{E(\dot S\dot X)\right\}^2\\
&-p_1E(\dot S^2)\left\{s_1\sum_j\mu_j^{-1}E(\dot X\dot X_j)-\sum_j\mu_j^{-2}E(\dot X\dot X_j)\right\}+O(n^{-1}).
\end{align*}
This, together with (\ref{d321}), (\ref{Ed02Hd11}), (\ref{Ed03d11}), (\ref{Ed002rmum1d11}), (\ref{md10}) and simplifications yields
\begin{align}
&E(d_0^2d_{321}d_{11})\notag\\
&=p_1E(d_0^2Hd_{11})+p_1s_1
\left[2E(d_0^2d_{00}d_{11})+E\{d_0^2(r,\mu^{m-1})d_{11}\}\right]
\notag\\
&+E(d_0^2C_{11}d_{11})\notag\\
&=E(\dot S^2)\left\{C_{10}\sum_j\mu_j^{-1}E(\tilde Z\dot X_j)-E(\tilde Z\dot X)\right\}-2E(\dot S\tilde Z)E(\dot S\dot X)\notag\\
&-p_1^2E(\dot S^2)\left\{s_1\sum_j\mu_j^{-1}E(\tilde Z\dot X_j)-\sum_j\mu_j^{-2}E(\tilde Z\dot X_j)\right\}\notag\\
&+s_1E(\dot S^2)\left\{
-6E(\dot S\dot X)-d_{10}p_1^2s_1-E(V\dot X)
-2p_1^3s_1^2+p_1^3s_2\right\}\notag\\
&+O(n^{-1})\notag\\
&=-2E(\dot S\tilde Z)E(\dot S\dot X)+E(\dot S^2)\left\{(-p_1d_{10}-2p_1^2s_1)\sum_j\mu_j^{-1}E(\tilde Z\dot X_j)
\right.\notag\\
&-E(\tilde Z \dot X)+p_1^2\sum_j\mu_j^{-2}E(\tilde Z\dot X_j)
-6s_1E(\dot S\dot X)\notag\\
&\left.-d_{10}p_1^2s_1^2-s_1E(V\dot X)
-2p_1^3s_1^3+p_1^3s_1s_2\right\}+O(n^{-1}).
\label{Ed02d321d11}
\end{align}

Hence, using (\ref{G1}) and (\ref{Ed02d321r1}),
\begin{align}
&E(d_0^2d_{321}G_1)\notag\\
&=\sum_j\mu_j^{-1}E(d_0^2d_{321}r_j)+d_{10}^{-1}E(d_0^2d_{321}d_{11})\notag\\
&=E(\dot S^2)\left\{p_1\sum_j\mu_j^{-1}E(\tilde Z\dot X_j)+p_1^2s_1^2\right\}+d_{10}^{-1}E(d_0^2d_{321}d_{11})\notag\\
&+O(n^{-1})\notag\\
&=-2d_{10}^{-1}E(\dot S\tilde Z)E(\dot S\dot X)+d_{10}^{-1}E(\dot S^2)
\left\{-2p_1^2s_1\sum_j\mu_j^{-1}E(\tilde Z\dot X_j)-E(\tilde Z\dot X)\right.\notag\\
&\left.+p_1^2\sum_j\mu_j^{-2}E(\tilde Z\dot X_j)-6s_1E(\dot S\dot X)
-s_1E(V\dot X)-2p_1^3s_1^3+p_1^3s_1s_2\right\}
\notag\\&+O(n^{-1}).\label{Ed02d321G1}
\end{align}

The next term needed is $E(d_0^2G_2)$. To this end, we have in the usual manner that
\begin{equation}
E(d_{00}^2r_1r_2)=n^{-2}\sum_i\sum_j\sum_k\sum_lE(U_iU_j\dot y_{k1}\dot y_{l2})=O(n^{-1}).\label{Ed002r1r2}
\end{equation}
Moreover, from (\ref{d11}), (\ref{Ed002r1r1}), (\ref{Ed002C11r1}) and (\ref{Ed002r1r2}),
\begin{align}
&E(d_0^2d_{11}r_1)\notag\\
&=p_1^{-1}\left\{C_{10}\mu_1^{-1}E(d_0^2r_1^2)-E(d_0^2C_{11}r_1)\right\}\notag\\
&-p_1\left\{s_1\mu_1^{-1}E(d_0^2r_1^2)-\mu_1^{-2}E(d_0^2r_1^2)\right\}+O(n^{-1})\notag\\
&=E(\dot S^2)\left[p_1^{-1}\left\{C_{10}-E(\dot X\dot X_1)\right\}-p_1\left(s_1-\mu_1^{-1}\right)\right]+O(n^{-1}).
\label{Ed02d11r1}
\end{align}
Now, via (\ref{G2}) and (\ref{Ed002r1r2}),
\begin{align*}
E(d_0^2G_2)&=3\sum_j\mu_j^{-2}E(d_{00}^2r_j^2)+4d_{10}^{-1}\sum_j\mu_j^{-1}E(d_{00}^2d_{11}r_j)\\
&+3d_{10}^{-2}E(d_{00}^2d_{11}^2)-2d_{10}^{-1}E(d_{00}^2d_{12})
+O(n^{-1}),
\end{align*}
and from (\ref{Ed002r1r1}), (\ref{Ed02d11r1}), (\ref{Ed002d112}), (\ref{Ed002d12}) and simplifications, we get
\begin{align}
E(d_0^2G_2)
&=3E(\dot S^2)s_1+d_{10}^{-1}E(\dot S^2)
\left\{6p_1^{-1}s_1C_{10}-6p_1^{-1}A_1
\right.\notag\\
&\left.-4p_1(s_1^2-s_2)\right\}+3d_{10}^{-2}E(d_{00}^2d_{11}^2)+O(n^{-1}).
\label{Ed02G2}
\end{align}
Moreover, (\ref{D}) and (\ref{Ed002d112}) yield
\begin{align*}
&d_{10}^{-2}E(d_{00}^2d_{11}^2)\\
&=2d_{10}^{-2}p_1^{-2}\left\{E(\dot S\dot X)\right\}^2+E(\dot S^2)
\left\{d_{10}^{-2}p_1^{-2}E(\dot X^2)+2d_{10}^{-1}p_1^{-1}A_1\right.\\
&+4d_{10}^{-2}s_1A_1+s_1+4d_{10}^{-1}p_1s_1^2+4d_{10}^{-2}p_1^2s_1^3
-2d_{10}^{-2}A_2+d_{10}^{-2}p_1^2s_3\\
&\left.-2d_{10}^{-1}p_1s_2
-4d_{10}^{-2}p_1^2s_1s_2
\right\}+O(n^{-1/2}),
\end{align*}
and together with (\ref{Ed02G2}), (\ref{md10}) and simplifications, we find
\begin{align}
&E(d_0^2G_2)\notag\\
&=6d_{10}^{-2}p_1^{-2}\left\{E(\dot S\dot X)\right\}^2+E(\dot S^2)
\left\{6s_1+6d_{10}^{-1}p_1^{-1}s_1C_{10}+8d_{10}^{-1}p_1s_1^2\right.\notag\\
&-2d_{10}^{-1}p_1s_2+3d_{10}^{-2}p_1^{-2}
E(\dot X^2)+12d_{10}^{-2}s_1A_1+12d_{10}^{-2}p_1^2s_1^3
\notag\\
&\left.-6d_{10}^{-2}A_2+3d_{10}^{-2}p_1^2s_3
-12d_{10}^{-2}p_1^2s_1s_2\right\}+O(n^{-1/2})\notag\\
&=6d_{10}^{-2}p_1^{-2}\left\{E(\dot S\dot X)\right\}^2+E(\dot S^2)
\left\{2d_{10}^{-1}p_1(s_1^2-s_2)\right.\notag\\
&+3d_{10}^{-2}p_1^{-2}
E(\dot X^2)+12d_{10}^{-2}s_1A_1+12d_{10}^{-2}p_1^2s_1^3
\notag\\
&\left.-6d_{10}^{-2}A_2+3d_{10}^{-2}p_1^2s_3
-12d_{10}^{-2}p_1^2s_1s_2\right\}+O(n^{-1/2}).
\label{Ed02G2b}
\end{align}
Now, by (\ref{s022}), (\ref{Ed02d332}), (\ref{Ed02d321G1}), (\ref{d320b}) and (\ref{Ed02G2b}),
\begin{align*}
&E(d_0^2s_{022})\\
&=E(d_0^2d_{322})-2E(d_0^2d_{321}G_1)+d_{320}E(d_0^2G_2)\\
&=E(\dot S^2)\{6E(\dot S Z)+E(ZV)\}+4d_{10}^{-1}E(\dot S\tilde Z)E(\dot S\dot X)\\
&-2d_{10}^{-1}E(\dot S^2)\left\{-2p_1^2s_1\sum_j\mu_j^{-1}E(\tilde Z\dot X_j)-E(\tilde Z\dot X)\right.\\
&\left.+p_1^2\sum_j\mu_j^{-2}E(\tilde Z\dot X_j)-6s_1E(\dot S\dot X)
-s_1E(V\dot X)-2p_1^3s_1^3+p_1^3s_1s_2\right\}\\
&-6d_{10}^{-1}p_1^{-1}\left\{E(\dot S\dot X)\right\}^2\\
&-d_{10}^{-1}E(\dot S^2)
\left\{2d_{10}p_1^2(s_1^2-s_2)+3p_1^{-1}
E(\dot X^2)+12p_1s_1A_1+12p_1^3s_1^3
\right.\\
&\left.-6p_1A_2+3p_1^3s_3
-12p_1^3s_1s_2\right\}+O(n^{-1/2}),
\end{align*}
and using $\tilde Z=Z+p_1^{-1}\dot X$, defining for $k=1,2$,
\begin{equation}
\tilde A_k=\sum_j\mu_j^{-k}E(\tilde Z\dot X_j)-3p_1^{-1}A_k
=\sum_j\mu_j^{-k}E\{(Z-2p_1^{-1}\dot X)\dot X_j\},
\label{tAk}
\end{equation}
we get after some simplification
\begin{align}
E(d_0^2s_{022})
&=E(\dot S^2)\{6E(\dot S Z)+E(ZV)\}+4d_{10}^{-1}E(\dot S Z)E(\dot S\dot X)\notag\\
&-2d_{10}^{-1}p_1^{-1}\left\{E(\dot S\dot X)\right\}^2-d_{10}^{-1}E(\dot S^2)\tilde F+O(n^{-1/2}),
\label{Ed02s022}
\end{align}
where 
\begin{align}
\tilde F&=-4p_1^2s_1\tilde A_1-2E(Z\dot X)
+p_1^{-1}E(\dot X^2)+2p_1^2\tilde A_2\notag\\
&-12s_1E(\dot S\dot X)-2s_1E(V\dot X)+2d_{10}p_1^2(s_1^2-s_2)
\notag\\
&+p_1^3(8s_1^3-10s_1s_2+3s_3).\label{tildeF}
\end{align}
To get a final expression for $E(f_0)=E(f_{01})+E(f_{02})$ (see (\ref{f0dec})), then in view of (\ref{Ef01}), (\ref{f02c}) and (\ref{d320b}),
we need
$$\left(p_1^{-1}d_{10}^{-1}+\frac{1}{2}p_1^{-2}d_{10}^{-2}d_{320}\right)E(d_0^2)=
\frac{1}{2}p_1^{-1}d_{10}^{-1}E(d_0^2),$$
and via (\ref{Ed02b}), this is
\begin{align}
&\left(p_1^{-1}d_{10}^{-1}+\frac{1}{2}p_1^{-2}d_{10}^{-2}d_{320}\right)E(d_0^2)\notag\\
&=\frac{1}{2}-\frac{1}{2}n^{-1}d_{10}^{-1}\sum_{j_1<...<j_{m-2}}\mu_{j_1}\cdots\mu_{j_{m-2}}+O(n^{-2}).\label{f0term1}
\end{align}
Moreover, via (\ref{Ed02s011}), (\ref{Ed02s021}) and simplifications,
\begin{align}
&n^{1/2}p_1^{-2}d_{10}^{-2}\left\{p_1d_{10}E(d_0^2s_{011})+\frac{1}{2}E(d_0^2s_{021})\right\}\notag\\
&=p_1^{-2}d_{10}^{-2}\left\{\frac{1}{2}E(\dot S^2\dot X)+\frac{1}{2}p_1E(\dot S^2 Z)+p_1^2s_1d_{10}+\frac{1}{2}p_1s_1\eta\right\}\notag\\
&+O(n^{-1/2}).\label{f0term2}
\end{align}
Furthermore, via (\ref{Ed002s012}), (\ref{Ed02s022}) and a little simplification, we get
\begin{align}
&p_1^{-2}d_{10}^{-2}\left\{p_1d_{10}E(d_0^2s_{012})+\frac{1}{2}E(d_0^2s_{022})\right\}\notag\\
&=p_1^{-2}d_{10}^{-2}\left[p_1^{-1}d_{10}^{-1}
\left\{E(\dot S\dot X)\right\}^2+3E(\dot S^2)E(\dot S Z)
+\frac{1}{2}E(\dot S^2)E(ZV)\right.\notag\\
&\left.+2d_{10}^{-1}E(\dot S Z)E(\dot S\dot X)+d_{10}^{-1}E(\dot S^2)
\left\{p_1d_{10}^2F-\frac{1}{2}\tilde F\right\}\right]+O(n^{-1/2}),\label{f0term3}
\end{align}
where via (\ref{F}), (\ref{tildeF}) and simplification,
\begin{align}
&p_1d_{10}^2F-\frac{1}{2}\tilde F\notag\\
&=\frac{1}{2}p_1^{-1}E(\dot X^2)+E(Z\dot X)+6s_1E(\dot S\dot X)+s_1E(V\dot X)\notag\\
&+2p_1^2s_1\sum_j\mu_j^{-1}E(Z\dot X_j)-
p_1^2\sum_j\mu_j^{-2}E(Z\dot X_j)\notag\\
&+p_1^3s_1s_2-\frac{1}{2}p_1^3s_3.
\label{Fdiff}
\end{align}
Hence, we may conclude from (\ref{f0dec}), (\ref{Ef01}), (\ref{f02c}) and (\ref{f0term1})-(\ref{Fdiff}) that
\begin{equation}
E(f_0)=\frac{1}{2}+n^{-1}\rho_0+O(n^{-3/2}),\label{Ef0}
\end{equation}
where
\begin{align}
\rho_0&=-\frac{1}{2}d_{10}^{-1}\sum_{j_1<...<j_{m-2}}\mu_{j_1}\cdots\mu_{j_{m-2}}+p_1^{-2}d_{10}^{-2}\cdot\notag\\
&\left[\frac{1}{2}E(\dot S^2\dot X)+\frac{1}{2}p_1E(\dot S^2 Z)+p_1^2s_1d_{10}+\frac{1}{2}p_1s_1\eta\right.\notag\\
&+p_1^{-1}d_{10}^{-1}
\left\{E(\dot S\dot X)\right\}^2+3E(\dot S^2)E(\dot S Z)
+\frac{1}{2}E(\dot S^2)E(ZV)\notag\\
&+2d_{10}^{-1}E(\dot S Z)E(\dot S\dot X)+d_{10}^{-1}E(\dot S^2)\cdot
\notag\\
&\left\{\frac{1}{2}p_1^{-1}E(\dot X^2)+E(Z\dot X)+6s_1E(\dot S\dot X)+s_1E(V\dot X)\right.\notag\\
&+2p_1^2s_1\sum_j\mu_j^{-1}E(Z\dot X_j)-
p_1^2\sum_j\mu_j^{-2}E(Z\dot X_j)\notag\\
&\left.\left.+p_1^3s_1s_2-\frac{1}{2}p_1^3s_3\right\}\right].\label{rho0}
\end{align}
To get further, it is time to evaluate the expectation terms in (\ref{rho0}). As for $\dot S$, it turns out easiest to replace the formula in (\ref{dotS}) by the equivalent
\begin{align}
\dot S&=\prod_j X_j-\dot X_1\prod_{j\neq 1}\mu_j-...-\dot X_m\prod_{j\neq m}\mu_j-\prod_j\mu_j\notag\\
&=\prod_j X_j-p_1\left(1+\sum_j\mu_j^{-1}\dot X_j\right).\label{dotSb}
\end{align}
This gives
\begin{align}
E(\dot S^2)
&=\prod_j E(X_j^2)+p_1^2\left\{1+\sum_j\mu_j^{-2}E(\dot X_j^2)\right\}\notag\\
&-2p_1\left\{\prod_j E(X_j)+\sum_j\mu_j^{-1}E\left(\dot X_j\prod_k X_k\right)\right\}\notag\\
&=p_1\prod_j(1+\mu_j)-p_1^2(1+s_1)\notag\\
&=p_1d_{10},\label{EdotS2}
\end{align}
comparing to (\ref{d10}) and 
\begin{equation}
C_{10}=E(X)=p_1\left\{p_1-\prod_j(1+\mu_j)\right\},\label{C10b}
\end{equation}
in the last step.

Moreover, because $X$ in (\ref{X})
is derived from $C_1$ in (\ref{C1}), we may write
\begin{equation}
X=\prod_j X_j\left\{\prod_j(X_j-1)-\prod_j X_j\right\}
=\prod_j X_j(X_j-1)-\prod_j X_j^2,\label{Xb}
\end{equation}
To calculate $E(\dot X^2)=E(X^2)-C_{10}^2$ (recall that $C_{10}=E(X)$), we need
$$E(X^2)=\prod_jE\left\{X_j^2(X_j-1)^2\right\}
+\prod_jE(X_j^4)-2\prod_jE\left\{X_j^3(X_j-1)\right\},$$
where knowledge of noncentral Poisson moments and simplifications leads to
$$E(X^2)=p_1^2p_{241}+p_1p_{1761}-2p_1^2p_{452},$$
letting
\begin{equation}
p_{abc}\eqd\prod_j(a+b\mu_j+c\mu_j^2+d\mu_j^3),\label{pabc}
\end{equation}
for integers $a,b,c,d$, where $p_{a,b,c,0}=p_{abc}$, etcetera.. Via (\ref{C10b}), it follows that
\begin{align}
&E(\dot X^2)\notag\\
&=p_1\left\{p_1p_{241}
+p_{1761}-2p_1p_{452}-p_1(p_1-p_{11})^2\right\}.\label{EdotX2}
\end{align}
Next, we work out $E(\dot S\dot X)=E(\dot S X)$. To this end, (\ref{dotSb}) and (\ref{Xb}) yield
\begin{align*}
E(\dot S\dot X)&=\prod_j E\left\{X_j^2(X_j-1)\right\}-\prod_j E(X_j^3)\\
&-p_1\left[\prod_j E\left\{X_j(X_j-1)\right\}+\sum_j\mu_j^{-1}
E\left\{\dot X_j\prod_k X_k(X_k-1)\right\}\right]\\
&+p_1\left\{\prod_j E(X_j^2)+\sum_j\mu_j^{-1}
E\left(\dot X_j\prod_k X_k^2\right)\right\},
\end{align*}
and it follows that
\begin{align}
&E(\dot S\dot X)\notag\\
&=p_1\{p_1p_{21}
-p_{131}-p_1^2-2p_1^2s_1
\notag\\
&+p_1p_{11}+p_1\sum_j\mu_j^{-1}(1+2\mu_j)\prod_{k\neq j}(1+\mu_k)\}\notag\\
&=p_1(p_1p_{21}
-p_{131}-p_1^2-2p_1^2s_1
+p_1p_{11}+p_1p_{11}s_{1;11}^{(12)}),\label{ESX}
\end{align}
where
\begin{equation}
s_{a;bc}^{(de)}\eqd\sum_j\mu_j^{-a}(b+c\mu_j)^{-1}(d+e\mu_j),\label{sabcd}
\end{equation}
for integers $a,b,c,d,e$, etcetera.
Next, we use (\ref{dotSb}) to write
\begin{align}
\dot S^2
&=\prod_j X_j^2-2p_1\left(\prod_j X_j+\sum_j\mu_j^{-1}\dot X_j\prod_k X_k\right)\notag\\
&+p_1^2\left(1+2\sum_j\mu_j^{-1}\dot X_j+\sum_j\sum_k\mu_j^{-1}\mu_k^{-1}\dot X_j\dot X_k\right),\label{dotS2}
\end{align}
which via (\ref{Xb}) yields
\begin{align*}
&E(\dot S^2 X)\\
&=\prod_j E\{X_j^3(X_j-1)\}-\prod_j E(X_j^4)\\
&-2p_1\left[\prod_j E\{X_j^2(X_j-1)\}-\prod_j E(X_j^3)\right]\\
&-2p_1\left[\sum_j\mu_j^{-1} E\left\{\dot X_j\prod_k X_k^2(X_k-1)\right\}
-\sum_j\mu_j^{-1} E\left(\dot X_j\prod_k X_k^3\right)\right]\\
&+p_1^2\left[\prod_jE\{X_j(X_j-1)\}-\prod_jE(X_j^2)\right]\\
&+2p_1^2\left[\sum_j\mu_j^{-1}E\left\{\dot X_j\prod_k X_k(X_k-1)\right\}
-\sum_j\mu_j^{-1}E\left(\dot X_j\prod X_k^2\right)\right]\\
&+p_1^2\left[\sum_j\sum_k\mu_j^{-1}\mu_k^{-1}
E\left\{\dot X_j\dot X_k\prod_l X_l(X_l-1)\right\}\right.\\
&\left.-\sum_j\sum_k\mu_j^{-1}\mu_k^{-1}E\left(\dot X_j\dot X_k\prod_l X_l^2\right)\right],
\end{align*}
and inserting the required moment formulae, we obtain
\begin{align}
&E(\dot S^2 X)\notag\\
&=p_1(p_1p_{451}
-p_{1761})-2p_1^2(p_1p_{21}-p_{131})\notag\\
&-2p_1^2\left\{p_1\sum_j\mu_j^{-1}(4+3\mu_j)\prod_{k\neq j}(2+\mu_k)\right.\notag\\
&\left.-\sum_j\mu_j^{-1}(1+6\mu_j+3\mu_j^2)\prod_{k\neq j}(1+3\mu_k+\mu_k^2)\right\}
+p_1^3(p_1-p_{11})\notag\\
&+2p_1^3\left\{2p_1s_1-\sum_j\mu_j^{-1}(1+2\mu_j)\prod_{k\neq j}(1+\mu_k)\right\}\notag\\
&+p_1^3\left\{p_1s_1(4s_1+1)
-\sum_j\mu_j^{-2}(1+5\mu_j+\mu_j^2)\prod_{k\neq j}(1+\mu_k)\right.\notag\\
&\left.-\sum_{j\neq k}\mu_j^{-1}\mu_k^{-1}(1+2\mu_j)(1+2\mu_k)\prod_{l\notin(j,k)}(1+\mu_l)\right\}\notag\\
&=p_1(p_1p_{451}
-p_{1761})-2p_1^2(p_1p_{21}-p_{131})-2p_1^2(p_1p_{21}s_{1;21}^{(43)}-p_{131}s_{1;131}^{(163)})\notag\\
&+p_1^3(p_1-p_{11})+2p_1^3(2p_1s_1-p_{11}s_{1;11}^{(12)})\notag\\
&+p_1^3\{p_1s_1(4s_1+1)-p_{11}s_{2;11}^{(151)}
-p_{11}(s_{1;11}^{(12)})^2+p_{11}s_{2;121}^{(144)}\}
.\label{ES2X}
\end{align}
Moreover, by (\ref{EdotS2}) and (\ref{C10b}),
\begin{align}
E(\dot S^2\dot X)&=E(\dot S^2 X)-C_{10}E(\dot S^2)\notag\\
&=E(\dot S^2 X)-p_1^2d_{10}\left\{p_1-\prod_j(1+\mu_j)\right\}\notag\\
&=E(\dot S^2 X)-p_1^2d_{10}(p_1-p_{11}).
\label{ESXdot}
\end{align}
Next, by (\ref{V}) we may write
\begin{equation}
V=p_1\sum_j\mu_j^{-1}\dot X_j,\label{Vb}
\end{equation}
and via (\ref{Xb}), we find
\begin{align}
&E(V\dot X)=E(VX)\notag\\
&=p_1\left[\sum_j\mu_j^{-1}
E\left\{\dot X_j\prod_k X_k(X_k-1)\right\}-\sum_j\mu_j^{-1}
E\left(\dot X_j\prod_k X_k^2\right)\right]\notag\\
&=p_1^2\left\{2p_1s_1-\sum_j\mu_j^{-1}(1+2\mu_j)
\prod_{k\neq j}(1+\mu_k)\right\}\notag\\
&=p_1^2(2p_1s_1-p_{11}s_{1;11}^{(12)})
.\label{EVX}
\end{align}
We are left with the terms involving $Z$. To attack these terms, observe that (\ref{Z}) may be reformulated as
\begin{equation}
Z=\sum_j\prod_{k\neq j}X_k-p_1s_1.\label{Zb}
\end{equation}
From (\ref{Zb}), we find
$$E(Z\dot X_1)=\sum_j\prod_{k\neq j}E(\dot X_1 X_k)
=\sum_{j\neq 1}\prod_{k\neq j}\mu_k=\sum_{j\neq 1}\mu_j^{-1}p_1=p_1(s_1-\mu_1^{-1}),
$$
which implies
\begin{equation}
\sum_j\mu_j^{-1}E(Z\dot X_j)=p_1(s_1^2-s_2),\label{ZXsum1}
\end{equation}
and
\begin{equation}
\sum_j\mu_j^{-2}E(Z\dot X_j)=p_1(s_1s_2-s_3),\label{ZXsum2}
\end{equation}
Moreover, (\ref{Vb}) and (\ref{Zb}) yield
\begin{align}
&E(ZV)\notag\\
&=p_1\sum_j\mu_j^{-1}\sum_kE\left(\dot X_j\prod_{l\neq k}X_l\right)=p_1\sum_j\mu_j^{-1}\sum_{k\neq j}E\left(\dot X_j\prod_{l\neq k}X_l\right)\notag\\
&=p_1\sum_j\mu_j^{-1}\sum_{k\neq j}\prod_{l\neq k}\mu_l
=p_1^2\sum_j\mu_j^{-1}\sum_{k\neq j}\mu_k^{-1}
=p_1^2(s_1^2-s_2).\label{EZV}
\end{align}
Furthermore, we get via (\ref{Xb}), (\ref{C10b}) and (\ref{Zb}) that
\begin{align}
&E(Z\dot X)=E(ZX)\notag\\
&=\sum_j E\left\{\prod_{k\neq j}X_k\prod_l X_l(X_l-1)\right\}
-\sum_j E\left(\prod_{k\neq j}X_k\prod_l X_l^2\right)-p_1s_1E(X)
\notag\\
&=\sum_j E\left\{X_j(X_j-1)\right\}\prod_{k\neq j}E\left\{X_k^2(X_k-1)\right\}\notag\\
&-\sum_j E(X_j^2)\prod_{k\neq j}E(X_k^3)-p_1s_1E(X)\notag\\
&=p_1^2\sum_j\prod_{k\neq j}(2+\mu_k)
-p_1\sum_j(1+\mu_j)\prod_{k\neq j}(1+3\mu_k+\mu_k^2)\notag\\
&-p_1^2s_1(p_1-p_{11})\notag\\
&=p_1^2p_{21}s_{21}-p_1p_{131}s_{0;131}^{(11)}
-p_1^2s_1(p_1-p_{11})
.\label{EZX}
\end{align}
Next, we may observe that from (\ref{dotSb}) and (\ref{Vb}), we have
\begin{equation}
\dot S=\prod_j X_j-p_1-V,\label{dotSc}
\end{equation}
and so, because $E(Z)=0$,
\begin{equation}
E(\dot SZ)=E\left(Z\prod_j X_j\right)-E(ZV),\label{ESZ}
\end{equation}
where via (\ref{Zb}),
\begin{align}
E\left(Z\prod_j X_j\right)&=\sum_j E\left(\prod_l X_l\prod_{k\neq j}X_k\right)-p_1s_1E\left(\prod_j X_j\right)\notag\\
&=\sum_j E(X_j)\prod_{k\neq j}E(X_k^2)-p_1^2s_1\notag\\
&=p_1\sum_j\prod_{k\neq j}(1+\mu_k)-p_1^2s_1
=p_1p_{11}s_{0;11}-p_1^2s_1.
\label{EZprodX}
\end{align}
Thus, (\ref{EZV}), (\ref{ESZ}) and (\ref{EZprodX}) imply
\begin{equation}
E(\dot S Z)=p_1p_{11}s_{0;11}-p_1^2(s_1+s_1^2-s_2).\label{ESZb}
\end{equation}
Next, (\ref{dotSc}) yields, because $E(Z)=0$,
\begin{align}
E(\dot S^2 Z)&=E\left(Z\prod X_j^2\right)-2p_1E\left(Z\prod X_j\right)-2E\left(ZV\prod X_j\right)\notag\\
&+2p_1E(ZV)+E(ZV^2).\label{ES2Z}
\end{align}
Here, by (\ref{Zb}),
\begin{align}
&E\left(Z\prod_j X_j^2\right)\notag\\
&=\sum_j E\left(\prod_l X_l^2\prod_{k\neq j}X_k\right)-p_1s_1E\left(\prod_j X_j^2\right)\notag\\
&=\sum_j E(X_j^2)\prod_{k\neq j}E(X_k^3)-p_1^2s_1\prod_j(1+\mu_j)\notag\\
&=p_1\sum_j(1+\mu_j)\prod_{k\neq j}(1+3\mu_k+\mu_k^2)
-p_1^2s_1\prod_j(1+\mu_j)\notag\\
&=p_1p_{131}s_{0;131}^{(11)}-p_1^2p_{11}s_1.\label{EZprodX2}
\end{align}
Moreover, via (\ref{Vb}) and (\ref{Zb}),
\begin{align}
&E\left(ZV\prod_j X_j\right)\notag\\
&=\sum_j E\left(V\prod_q X_q\prod_{k\neq j}X_k\right)-p_1s_1E\left(V\prod_j X_j\right)\notag\\
&=p_1\left\{\sum_j\sum_l\mu_l^{-1}E\left(\dot X_l\prod_q X_q\prod_{k\neq j}X_k\right)-p_1s_1\sum_k\mu_k^{-1}E\left(\dot X_k\prod_j X_j\right)\right\}\notag\\
&=p_1\left\{\sum_j\mu_j^{-1}E(\dot X_j X_j)\prod_{k\neq j}E(X_k^2)+\sum_{j\neq l}\mu_l^{-1}E(X_j)\prod_{k\neq j}E(\dot X_k X_k^2)\right.\notag\\
&\left.-p_1s_1\sum_k\mu_k^{-1}E(\dot X_k X_k)\prod_{j\neq k}E(X_j)\right\}\notag\\
&=p_1^2\left\{\sum_j\mu_j^{-1}\prod_{k\neq j}(1+\mu_k)
+\sum_{j\neq l}\mu_l^{-1}\prod_{k\neq j}(1+2\mu_k)
-p_1s_1^2\right\}\notag\\
&=p_1^2\{p_{11}s_{1;11}+p_{12}(s_1s_{12}-s_{1;12})-p_1s_1^2\}
.\label{EZVprodX}
\end{align}
Furthermore, (\ref{Zb}) gives
\begin{equation}
E(ZV^2)=\sum_j E\left(V^2\prod_{k\neq j}X_k\right)-p_1s_1E(V^2),\label{EZV2}
\end{equation}
where from (\ref{Vb}),
\begin{align}
&E(V^2) \notag\\
&=p_1^2\sum_j\sum_k\mu_j^{-1}\mu_k^{-1}E(\dot X_j\dot X_k)=p_1^2\sum_j\mu_j^{-2}E(\dot X_j^2)=p_1^2 s_1,
\label{EV2}
\end{align}
and
\begin{align}
\sum_j E\left(V^2\prod_{k\neq j}X_k\right)&=
p_1^2\sum_j\sum_q\sum_l\mu_q^{-1}\mu_l^{-1}E\left(\dot X_q\dot X_l\prod_{k\neq j}X_k\right)\notag\\
&=p_1^2(v_1+v_2),\label{EV2term1}
\end{align}
where we get the term with $q=l$ as
\begin{align}
v_1&=\sum_j\sum_l\mu_l^{-2}E\left(\dot X_l^2\prod_{k\neq j}X_k\right)
\notag\\
&=\sum_j\mu_j^{-2}E(\dot X_j^2)\prod_{k\neq j}E(X_k)
+\sum_{j\neq l}\mu_l^{-2}E(\dot X_l^2 X_l)\prod_{k\notin(j,l)}E(X_k)
\notag\\
&=p_1s_2+p_1\sum_{j\neq l}\mu_l^{-2}\mu_j^{-1}(1+\mu_l)=p_1s_2+p_1(s_1s_2^{(11)}-s_3^{(11)}),\label{v1}
\end{align}
and for $q\neq l$ in (\ref{EV2term1}) we have
\begin{align}
v_2&=\sum_j\sum_{q\neq l}\mu_q^{-1}\mu_l^{-1}E\left(\dot X_q\dot X_l\prod_{k\neq j}X_k\right)\notag\\
&=\sum_j\sum_{q\notin(j,l),j\neq l}\mu_q^{-1}\mu_l^{-1}E(\dot X_q X_q)E(\dot X_l X_l)\prod_{k\notin(q,l)}E(X_k)\notag\\
&=p_1\sum_j\sum_{q\notin(j,l),j\neq l}\mu_q^{-1}\mu_l^{-1}
=p_1(m-2)(s_1^2-s_2).\label{v2}
\end{align}
Hence, (\ref{EZV2})-(\ref{v2}) together yield
\begin{equation}
E(ZV^2)=p_1^3\left\{(m-3)(s_1^2-s_2)+s_1s_2^{(11)}-s_3^{(11)}\right\}.\label{EZV2b}
\end{equation}
Thus, from (\ref{ES2Z}), (\ref{EZprodX2}), (\ref{EZprodX}), (\ref{EZVprodX}), (\ref{EZV}), (\ref{EZV2b}) and some simplification,
\begin{align}
&E(\dot S^2Z)\notag\\
&=p_1p_{131}s_{0;131}^{(11)}-p_1^2\{p_{11}(s_1+2s_{0;11}+2s_{1;11})+2p_{12}(s_1s_{12}-s_{1;12})\}\notag\\
&+p_1^3\{2s_1+2s_1^2+(m-1)(s_1^2-s_2)+s_1s_2^{(11)}-s_3^{(11)}\}.\label{ES2Zb}
\end{align}
We may now insert (\ref{EdotS2}), (\ref{ZXsum1}), (\ref{ZXsum2}) and (\ref{EZV}) into (\ref{rho0}) to obtain, after small simplifications,
\begin{align}
\rho_0&=-\frac{1}{2}d_{10}^{-1}\sum_{j_1<...<j_{m-2}}\mu_{j_1}\cdots\mu_{j_{m-2}}+p_1^{-2}d_{10}^{-2}\cdot\notag\\
&\left[\frac{1}{2}E(\dot S^2\dot X)+\frac{1}{2}p_1E(\dot S^2 Z)+p_1^2s_1d_{10}+\frac{1}{2}p_1s_1\eta\right.\notag\\
&+p_1^{-1}d_{10}^{-1}
\left\{E(\dot S\dot X)\right\}^2+3p_1d_{10}E(\dot S Z)
+\frac{1}{2}p_1^3d_{10}(s_1^2-s_2)\notag\\
&+2d_{10}^{-1}E(\dot S Z)E(\dot S\dot X)+p_1\cdot
\notag\\
&\left\{\frac{1}{2}p_1^{-1}E(\dot X^2)+E(Z\dot X)+6s_1E(\dot S\dot X)+s_1E(V\dot X)\right.\notag\\
&\left.\left.+2p_1^3s_1(s_1^2-s_2)+\frac{1}{2}p_1^3s_3\right\}\right],\label{rho0b}
\end{align}
where
\begin{equation}\sum_{j_1<...<j_{m-2}}\mu_{j_1}\cdots\mu_{j_{m-2}}
=p_1\sum_{j<k}\mu_j^{-1}\mu_k^{-1}=\frac{1}{2}p_1(s_1^2-s_2)\label{musum}
\end{equation}
and from (\ref{eta}),
\begin{equation}\eta=p_1^2(p_{11}s_{1;11}-p_1s_2).\label{etac}
\end{equation}

\newpage

Next, we focus on $f_1$. We have from (\ref{f1}) that
\begin{equation}
f_1d_0^{-3}=\frac{d_2}{P^2d_1^4}D_1+D_2,\label{f1b}
\end{equation}
where
\begin{equation}
D_1=P(-d_1)+d_{32},\label{D1}
\end{equation}
and
\begin{equation}
D_2=\frac{d_{33}}{3P^3(-d_1)^3}.\label{D2}
\end{equation}
Now, from (\ref{P}), (\ref{md1}), (\ref{d32b}), (\ref{d320b}) and (\ref{D1}), we have
\begin{align}
&D_1\notag\\
&=p_1d_{10}\left(1+n^{-1/2}\sum_j\frac{r_j}{\mu_j}\right)
\left(1+n^{-1/2}\frac{d_{11}}{d_{10}}\right)+d_{320}+n^{-1/2}d_{321}
+O(n^{-1})\notag\\
&=n^{-1/2}\left(p_1d_{10}\sum_j\mu_j^{-1}r_j+p_1d_{11}+d_{321}\right)+O(n^{-1}),\label{D1b}
\end{align}
and with
\begin{equation}
d_{33}=d_{330}+n^{-1/2}d_{331}+O(n^{-1}),\label{d33b}
\end{equation}
we have by (\ref{D2}) that
\begin{align}
&D_2\notag\\
&=\frac{1}{3}p_1^{-3}d_{10}^{-3}d_{330}
\left(1+n^{-1/2}\sum_j\frac{r_j}{\mu_j}\right)^{-3}
\left(1+n^{-1/2}\frac{d_{11}}{d_{10}}\right)^{-3}\left(1+n^{-1/2}\frac{d_{331}}{d_{330}}\right)\notag\\&+O(n^{-1})\notag\\
&=\frac{1}{3}p_1^{-3}d_{10}^{-3}d_{330}
\left\{1+n^{-1/2}\left(-3\sum_j\mu_j^{-1}r_j-3d_{10}^{-1}d_{11}+
d_{330}^{-1}d_{331}\right)\right\}\notag\\&+O(n^{-1}).\label{D2b}
\end{align}
Since $D_1$ is of order $n^{-1/2}$, we only need to keep track of the main term (limit) of $d_2$. From (\ref{d2}), this limit is
\begin{equation}
d_{20}=p_1^{-2}C_{20}-p_1\sum_{j<k}\mu_j^{-1}\mu_k^{-1}
=p_1^{-2}C_{20}-\frac{1}{2}p_1(s_1^2-s_2),\label{d20}
\end{equation}
where $C_{20}$ is the limit of (\ref{C2}). In fact, (\ref{b2}) and (\ref{C2}) yield
\begin{align*}
C_{20}
&=E\left(\prod_j X_j\left[p_1s_1\left\{\prod_j(X_j-1)-\prod_j X_j\right\}\right.\right.\notag\\
&+\frac{1}{2}\prod_j(X_j-1)(X_j-2)-\frac{1}{2}\prod_jX_j(X_j-1)\notag\\
&\left.\left.+\prod_j X_j^2-\prod_j X_j(X_j-1)\right]\right)\notag\\
&=p_1s_1\prod_j E\{X_j(X_j-1)\}-p_1s_1\prod_j E(X_j^2)\notag\\
&+\frac{1}{2}\prod_j E\{X_j(X_j-1)(X_j-2)\}-\frac{3}{2}\prod_j E\{X_j^2(X_j-1)\}\notag\\
&+\prod_j E(X_j^3),
\end{align*}
It follows that
\begin{align*}
C_{20}&=p_1^3s_1-p_1^2s_1\prod_j(1+\mu_j)+\frac{1}{2}p_1^3\\
&-\frac{3}{2}p_1^2\prod_j(2+\mu_j)+p_1\prod_j(1+3\mu_j+\mu_j^2)\\
&=p_1^3s_1-p_1^2p_{11}s_1+\frac{1}{2}p_1^3-\frac{3}{2}p_1^2p_{21}+p_1p_{131},
\end{align*}
which via (\ref{d20}) yields
\begin{equation}
d_{20}=p_1s_1-p_{11}s_1+\frac{1}{2}p_1-\frac{3}{2}p_{21}+p_1^{-1}p_{131}-\frac{1}{2}p_1(s_1^2-s_2).\label{d20b}
\end{equation}

Now, by (\ref{f1b}) and (\ref{D1b}),
\begin{align}
f_1d_0^{-3}&=n^{-1/2}\frac{d_{20}}{p_1^2d_{10}^4}
\left(p_1d_{10}\sum_j\mu_j^{-1}r_j+p_1d_{11}+d_{321}\right)\notag\\
&+D_2+O(n^{-1}).\label{f1c}
\end{align}
To tackle $D_2$, we at first need $d_{330}$, the limit of $d_{33}$. Here, we find that from (\ref{C1}) and (\ref{e2i}), the limit of $\bar e_2$ is $C_{10}$, given in (\ref{C10b}). Analogously, define $e_{30}$ as the limit of $\bar e_3$, cf (\ref {e3i}). This yields
\begin{align}
e_{30}&=\frac{1}{2}\prod_j E\{X_j(X_j-1)(X_j-2)-\frac{3}{2}\prod_j E\{X_j^2(X_j-1)\}+\prod_j E(X_j^3)\notag\\
&=\frac{1}{2}p_1^3-\frac{3}{2}p_1^2\prod_j(2+\mu_j)+p_1\prod_j(1+3\mu_j+\mu_j^2)\notag\\
&=\frac{1}{2}p_1^3-\frac{3}{2}p_1^2p_{21}+p_1p_{131}.\label{e30}
\end{align}

Then, letting $M(2^{m-1},0^1)$ be the limit of $\by(2^{m-1},0^1)$ etcetera, (\ref{d33}) yields
\begin{align}
d_{330}&=2p_1M(2^{m-1},0^1)+3p_1M(2^{m-2},1^2)
\notag\\&+3M(1^{m-1},0^1)C_{10}+e_{30}.\label{d330}
\end{align}
Here,
\begin{align} 
M(2^{m-1},0^1)&=\sum_{j_1<...<j_{m-1}}\prod_k\mu_{j_k}^2=p_1^2s_2,\label{M2m101}\\
M(2^{m-2},1^2)&=p_1\sum_{j_1<...<j_{m-2}}\prod_k\mu_{j_k}
=p_1^2\sum_{j<k}\mu_j^{-1}\mu_k^{-1}\notag\\
&=\frac{1}{2}p_1^2(s_1^2-s_2),\label{M2m212}\\
M(1^{m-1},0^1)&=p_1s_1.\label{M1m101}
\end{align} 
Inserting this together with (\ref{C10b}) and (\ref{e30}) into (\ref{d330}), we get
\begin{equation}
d_{330}=p_1^3\left(\frac{3}{2}s_1^2+\frac{1}{2}s_2\right)
+3p_1^2s_1(p_1-p_{11})+\frac{1}{2}p_1^3-\frac{3}{2}p_1^2p_{21}+p_1p_{131}.\label{d330b}
\end{equation}

To find $d_{331}$, we need the expansions
\begin{equation}
\bar e_2=C_{10}+n^{-1/2}C_{11}+O(n^{-1})\label{bare2}
\end{equation}
and
\begin{equation}
\bar e_3=e_{30}+n^{-1/2}e_{31}+O(n^{-1}),\label{bare3}
\end{equation}
where $e_{30}$ is in (\ref{e30}).
Moreover, we write e.g.
$$\by(2^{m-1},0^1)=M(2^{m-1},0^1)+n^{-1/2}\br(2^{m-1},0^1).$$
Next, analogous to (\ref{d32term1}), we have
$$3\bar y_{11...1}-P=2p_1+n^{-1/2}\{3d_{00}+2(r,\mu^{m-1})\}+O(n^{-1}),$$
and similarly,
$$\bar y_{11...1}=p_1+n^{-1/2}\{d_{00}+(r,\mu^{m-1})\}+O(n^{-1}),$$
implying via (\ref{d33}) and (\ref{M2m101})-(\ref{M1m101}) that
\begin{align}
d_{331}&=M(2^{m-1},0^1)\{3d_{00}+2(r,\mu^{m-1})\}
+2p_1\br(2^{m-1},0^1)\notag\\
&+3M(2^{m-2},1^2)\{d_{00}+(r,\mu^{m-1})\}+3p_1\br(2^{m-2},1^2)\notag\\
&+3M(1^{m-1},0^1)C_{11}+3C_{10}R(1^{m-1},0^1)+e_{31}\notag\\
&=p_1^2s_2\{3d_{00}+2(r,\mu^{m-1})\}
+2p_1\br(2^{m-1},0^1)\notag\\
&+\frac{3}{2}p_1^2(s_1^2-s_2)\{d_{00}+(r,\mu^{m-1})\}+3p_1\br(2^{m-2},1^2)\notag\\
&+3p_1s_1C_{11}+3C_{10}R(1^{m-1},0^1)+e_{31}.
\label{d331}
\end{align}
To find $E(d_0^3D_2)$, we need to calculate
\begin{align}
E(d_0^3r_1)&=E(d_{00}^3r_1)+O(n^{-1/2})\notag\\
&=n^{-2}\sum_i\sum_j\sum_k\sum_l E(U_iU_jU_k\dot y_{l1})+O(n^{-1/2})
\notag\\
&=O(n^{-1/2}).\label{Ed03r1}
\end{align}
Moreover, by (\ref{Ed03d11}) and (\ref{d0b}),
\begin{equation}
E(d_0^3d_{11})=-3p_1^{-1}E(\dot S^2)E(\dot S\dot X)+O(n^{-1/2}).
\label{Ed03d11b}
\end{equation}
To get $E(d_0^3d_{321})$, we find in the usual manner
\begin{align}
E(d_0^3H)&=n^{-2}\sum_i\sum_j\sum_k\sum_l E(U_iU_jU_kH_l)+O(n^{-1/2})\notag\\
&=3E(\dot S^2)E(\dot S Z)+O(n^{-1/2}),\label{Ed03H}
\end{align}
\begin{align}
E(d_0^3d_{00})&=n^{-2}\sum_i\sum_j\sum_k\sum_l E(U_iU_jU_kU_l)+O(n^{-1/2})\notag\\
&=3\{E(\dot S^2)\}^2+O(n^{-1/2}),\label{Ed03d00}
\end{align}
and
\begin{align}
E\{d_0^3(r,\mu^{m-1}\}&=n^{-2}\sum_i\sum_j\sum_k\sum_l E(U_iU_jU_kV_l)+O(n^{-1/2})\notag\\
&=O(n^{-1/2}),\label{Ed03rmum1}
\end{align}
which together with (\ref{Ed003C11}) inserted into (\ref{d321}) yields
\begin{align}
&E(d_0^3d_{321})\notag\\
&=p_1E(d_0^3H)+p_1s_1
\left[2E(d_0^3d_{00})+E\{d_0^3(r,\mu^{m-1})\}\right]+E(d_0^3C_{11})
\notag\\
&=3E(\dot S^2)\left\{p_1E(\dot SZ)+2p_1s_1E(\dot S^2)+E(\dot S\dot X)\right\}+O(n^{-1/2}).\label{Ed03d321}
\end{align}
To get $E(d_0^3D_2)$, we in addition need $E(d_0^3d_{331})$, and to that end, we must calculate terms like $E\{d_0^3R(1^{m-1},0^1)\}$. However, because
$$\by(1^{m-1},0^1)=\bar y_{1...10}+...+\bar y_{01...1},$$
where e.g.
$$\bar y_{1...10}=\mu_1\cdots\mu_{m-1}+n^{-1/2}(\mu_1\cdots\mu_{m-2}r_{m-1}+...+\mu_2\cdots\mu_{m-1}r_1)+O(n^{-1}),$$
it follows that
\begin{align*}
R(1^{m-1},0^1)&=\mu_1\cdots\mu_{m-2}r_{m-1}+...+\mu_2\cdots\mu_{m-1}r_1+...\\
&+\mu_2\cdots\mu_{m-1}r_m+...+\mu_3\cdots\mu_m r_2+O(n^{-1}),
\end{align*}
and so, it follows from (\ref{Ed03r1}) that
$$E\{d_0^3R(1^{m-1},0^1)\}=O(n^{-1/2}).$$
It is similarly seen that $E\{d_0^3\br(2^{m-1},0^1)\}$ and $E\{d_0^3\br(2^{m-2},1^2)\}$ are also $O(n^{-1/2})$.

However, the term $E(d_0^3e_{31})$ is not that simple, but we can perform the same 'trick' as with the corresponding term for $e_2$, $E(d_0^3C_{11})$. Indeed, as in (\ref{Ed003C11}), we have
\begin{align}
E(d_0^3e_{31})&=n^{-2}\sum_i\sum_j\sum_k\sum_l E(U_iU_jU_ke_{3l})+O(n^{-1/2})\notag\\
&=3E(\dot S^2)E(\dot S\epsilon_3)+O(n^{-1/2}),\label{Ed03e31}
\end{align}
where via (\ref{e3i}),
\begin{equation}
\epsilon_3=\frac{1}{2}\prod_j X_{j(3)}-\frac{3}{2}\prod_j X_j\prod_j X_{j(2)}
+\left(\prod_j X_j\right)^3,\label{eps3}
\end{equation}
where as before, $X_j$ is a Poisson variate with parameter $\mu_j$ and all $X_j$ are simultaneously independent. Moreover, $X_{j(2)}=X_j(X_j-1)$ and
\\ $X_{j(3)}=X_j(X_j-1)(X_j-2)$.
Moreover, we find in the usual manner that
\begin{align}
&E(\dot S\epsilon_3)\notag\\
&=\frac{1}{2}\prod_j E\{X_j^2(X_j-1)(X_j-2)\}-\frac{3}{2}\prod_jE\{X_j^3(X_j-1)\}+\prod_j E(X_j^4)\notag\\
&-p_1\left[\frac{1}{2}\prod_j E\{X_j(X_j-1)(X_j-2)\}-\frac{3}{2}\prod_jE\{X_j^2(X_j-1)\}+\prod_j E(X_j^3)\right.\notag\\
&+\frac{1}{2}\sum_j\mu_j^{-1}E\left\{\dot X_j\prod_k X_k(X_k-1)(X_k-2)\right\}\notag\\
&\left.-\frac{3}{2}\sum_j\mu_j^{-1}E\left\{\dot X_j\prod_k X_k^2(X_k-1)\right\}
+\sum_j\mu_j^{-1}E\left(\dot X_j\prod_k X_k^3\right)\right]\notag\\
&=\frac{1}{2}p_1^3\prod_j(3+\mu_j)-\frac{3}{2}p_1^2\prod_j(4+5\mu_j+\mu_j^2)\notag\\
&+p_1\prod_j(1+7\mu_j+6\mu_j^2+\mu_j^3)-\frac{1}{2}p_1^4+\frac{3}{2}p_1^3\prod_j(2+\mu_j)\notag\\
&-p_1^2\prod_j(1+3\mu_j+\mu_j^2)-\frac{3}{2}p_1^4s_1+\frac{3}{2}p_1^3\sum_j\mu_j^{-1}(4+3\mu_j)\prod_{k\neq j}(2+\mu_k)\notag\\
&-p_1^2\sum_j\mu_j^{-1}(1+6\mu_j+3\mu_j^2)\prod_{k\neq j}(1+3\mu_k+\mu_k^2)\notag\\
&=\frac{1}{2}p_1^3p_{31}-\frac{3}{2}p_1^2p_{451}+p_1p_{1761}-\frac{1}{2}p_1^4+\frac{3}{2}p_1^3p_{21}-p_1^2p_{131}-\frac{3}{2}p_1^4s_1\notag\\
&+\frac{3}{2}p_1^3p_{21}s_{1;21}^{(43)}-p_1^2p_{131}s_{1;131}^{(163)}.\label{ESeps3}
\end{align}

Now, from (\ref{d331}),
\begin{align*}
E(d_0^3d_{331})&=p_1^2s_2\left[3E(d_0^3d_{00})+2E\{d_0^3(r,\mu^{m-1}\}\right]\\
&+\frac{3}{2}p_1^2(s_1^2-s_2)\left[E(d_0^3d_{00})+E\{d_0^3(r,\mu^{m-1}\}\right]\\
&+3p_1s_1E(d_0^3C_{11})+E(d_0^3e_{31})+O(n^{-1/2}),
\end{align*}
and (\ref{Ed03d00}), (\ref{Ed03rmum1}), (\ref{Ed003C11}) and (\ref{eps3}) yield
\begin{align}
E(d_0^3d_{331})&=3E(\dot S^2)\left\{3p_1^2s_2E(\dot S^2)\right.+\frac{3}{2}p_1^2(s_1^2-s_2)E(\dot S^2)\notag\\
&\left.+3p_1s_1E(\dot S \dot X)+E(\dot S\epsilon_3)\right\}
+O(n^{-1/2}).\label{Ed03d331}
\end{align}
Now, via (\ref{D2b}) and (\ref{Ed03r1}),
\begin{align*}
E(d_0^3D_2)&=\frac{1}{3}p_1^{-3}d_{10}^{-3}d_{330}
\left[E(d_0^3)\right.\\
&\left.+n^{-1/2}\left\{-3d_{10}^{-1}E(d_0^3d_{11})+d_{330}^{-1}E(d_0^3d_{331})\right\}\right]+O(n^{-1}),
\end{align*}
and inserting (\ref{Ed02d00b}), (\ref{Ed03d11b}), (\ref{Ed03d331}) and (\ref{EdotS2}), we find
\begin{align}
&E(d_0^3D_2)\notag\\
&=n^{-1/2}\frac{1}{3}p_1^{-3}d_{10}^{-3}
\left[\left\{p_1d_{10}+9E(\dot S \dot X)\right\}d_{330}\right.\notag\\
&+3p_1d_{10}\left\{\frac{3}{2}p_1^3d_{10}(s_1^2+s_2)\right.\notag\\
&\left.\left.+3p_1s_1E(\dot S \dot X)+E(\dot S\epsilon_3)\right\}\right]+O(n^{-1}),\label{Ed03D2}
\end{align}
where $d_{330}$ is given in (\ref{d330b}), $E(\dot S X)$ is given in (\ref{ESX}) and $E(\dot S\epsilon_3)$ is given in (\ref{ESeps3}).

Now, we have from (\ref{f1c}) and (\ref{Ed03r1}) that
$$
E(f_1)=n^{-1/2}\frac{d_{20}}{p_1^2d_{10}^4}
\left\{p_1E(d_0^3d_{11})+E(d_0^3d_{321})\right\}
+E(d_0^3D_2)+O(n^{-1}),
$$
and (\ref{Ed03d11b}), (\ref{Ed03d321}) and a little simplification yields
\begin{align*}
E(f_1)&=n^{-1/2}\frac{d_{20}}{p_1d_{10}^4}3E(\dot S^2)
\left\{E(\dot S Z)+2s_1E(\dot S^2)\right\}\\
&+E(d_0^3D_2)+O(n^{-1}),
\end{align*}
and by (\ref{EdotS2}) and (\ref{Ed03D2}), we obtain
\begin{align}
n^{1/2}E(f_1)&=3d_{20}d_{10}^{-3}\left\{E(\dot S Z)+2p_1d_{10}s_1\right\}\notag\\
&+\frac{1}{3}p_1^{-2}d_{10}^{-2}d_{330}+3p_1^{-3}d_{10}^{-3}E(\dot S\dot X)d_{330}\notag\\
&+\frac{3}{2}p_1d_{10}^{-1}(s_1^2+s_2)+3p_1^{-1}d_{10}^{-2}s_1E(\dot S\dot X)+p_1^{-2}d_{10}^{-2}E(\dot S\epsilon_3)\notag\\&+O(n^{-1/2}).\label{f1}
\end{align}

Finally, we evaluate $f_2$. To this end, only $E(d_0^4)$ and limits of the terms in (\ref{f2}) are needed. At first, in the usual manner we have
\begin{align}
E(d_0^4)&=n^{-2}\sum_i\sum_j\sum_k\sum_lE(U_iU_jU_kU_l)+O(n^{-1/2})\notag\\
&=3\{E(\dot S^2)\}^2+O(n^{-1/2}).\label{Ed04}
\end{align}
Hence, from (\ref{f2}), (\ref{d320b}) and (\ref{EdotS2}), letting $d_{30}$ be the limit of $d_3$,
\begin{align}
E(f_2)&=\frac{3\{E(\dot S^2)\}^2}{4p_1^4d_{10}^6}\notag\\
&(8p_1^3d_{10}d_{20}^2+10p_1^2d_{20}^2d_{320}+4p_1d_{330}d_{10}d_{20}+d_{340}d_{10}^2
\notag\\
&+4p_1^3d_{10}^2d_{30}+4p_1^2d_{10}d_{30}d_{320})+O(n^{-1/2})\notag\\
&=\frac{3}{4p_1^2d_{10}^3}(-2p_1^3d_{20}^2+4p_1d_{330}d_{20}+d_{340}d_{10})\notag\\
&+O(n^{-1/2}),
\label{Ef2b}
\end{align}
Here, $d_{20}$ is given by (\ref{d20b}), $d_{330}$ is in (\ref{d330b}), and similarly, $d_{340}$ is obtained as the limit of $d_{34}$ in (\ref{d34}) as
\begin{align}
d_{340}&=3p_1M(3^{m-1},0^1)
\notag\\
&+4p_1\left\{M(3^{m-2},2^1,1^1)+M(3^{m-3},2^3)\right\}\notag\\
&+2\left\{M(1^{m-1},0^1)^2+2M(2^{m-1},0^1)+2M(2^{m-2},1^2)\right\}C_{10}\notag\\
&+4M(1^{m-1},0^1)e_{30}+e_{40},\label{d340}
\end{align}
where $e_{40}$ is the limit of (\ref{e4i}). Here, we get in the usual manner
\begin{align}
M(3^{m-1},0^1)&=p_1^3s_3,\label{M3m101}\\
M(3^{m-2},2^1,1^1)&=p_1^3\sum_j\mu_j^{-1}\sum_{k\neq j}\mu_k^{-2}=p_1^3(s_1s_2-s_3),\label{M3m22111}\\
M(3^{m-3},2^3)&=p_1^3\sum_{j_1<j_2<j_3}\prod_k\mu_{j_k}^{-1}=\frac{1}{6}p_1^3(s_1^3-3s_1s_2+2s_3),\label{M3m323}
\end{align}
and inserting (\ref{M3m101})-(\ref{M3m323}) together with (\ref{M2m101})-(\ref{M1m101}) into (\ref{d340}) and simplifying, we find
\begin{align}
d_{340}
&=p_1^4\left(\frac{2}{3}s_1^3+2s_1s_2+\frac{1}{3}s_3\right)+2p_1^2(2s_1^2+s_2)C_{10}\notag\\
&+4p_1s_1e_{30}+e_{40}.\label{d340b}
\end{align}
Here, (\ref{e4i}) implies
\begin{align}
&e_{40}\notag\\
&=\frac{1}{6}\prod_j E\{X_j(X_j-1)(X_j-2)(X_j-3)\}-\frac{2}{3}\prod_j E\{X_j^2(X_j-1)(X_j-2)\}\notag\\
&-\frac{1}{2}\prod_j E\{X_j^2(X_j-1)^2\}+2\prod_j E\{X_j^3(X_j-1)\}-\prod_j E(X_j^4)\notag\\
&=\frac{1}{6}p_1^4-\frac{2}{3}p_1^3\prod_j(3+\mu_j)
-\frac{1}{2}p_1^2\prod_j(2+4\mu_j+\mu_j^2)\notag\\
&+2p_1^2\prod_j(4+5\mu_j+\mu_j^2)
-p_1\prod_j(1+7\mu_j+6\mu_j^2+\mu_j^3)\notag\\
&=\frac{1}{6}p_1^4-\frac{2}{3}p_1^3p_{31}-\frac{1}{2}p_1^2p_{241}
+2p_1^2p_{451}-p_1p_{1761}.\label{e40}
\end{align}
Hence, inserting (\ref{C10b}), (\ref{e30}) and (\ref{e40}) into (\ref{d340b}), we get
\begin{align}
&d_{340}\notag\\
&=p_1^4\left(\frac{2}{3}s_1^3+2s_1s_2+\frac{1}{3}s_3\right)+2p_1^3(2s_1^2+s_2)(p_1-p_{11})\notag\\
&+2p_1^4s_1-6p_1^3p_{21}s_1+4p_1^2p_{131}s_1+\frac{1}{6}p_1^4-\frac{2}{3}p_1^3p_{31}
-\frac{1}{2}p_1^2p_{241}\notag\\
&+2p_1^2p_{451}-p_1p_{1761}.\label{d340c}
\end{align}
Hence, by inserting (\ref{d20b}), (\ref{d330b}) and (\ref{d340c}) into (\ref{Ef2b}), we get $E(f_2)$. One final simplification is to write
\begin{align}
q_1&=p_1s_1-p_{11}s_1-\frac{1}{2}p_1(s_1^2-s_2),\label{Q1}\\
q_2&=\frac{1}{2}p_1-\frac{3}{2}p_{21}+p_1^{-1}p_{131},\label{Q2}\\
q_3&=p_1^3\left(\frac{3}{2}s_1^2+\frac{1}{2}s_2\right)+3p_1^2s_1(p_1-p_{11}),\label{Q3}
\end{align}
so that by (\ref{d20b}) and (\ref{d330b}),
\begin{align}
d_{20}&=q_1+q_2,\label{d20c}\\
d_{330}&=q_3+p_1^2q_2,\label{d330c}
\end{align}
which yields
\begin{align}
&E(f_2)\notag\\
&=\frac{3}{4p_1^2d_{10}^3}\notag\\
&\cdot\left\{-2p_1^3(q_1+q_2)^2+4p_1(q_3+p_1^2q_2)(q_1+q_2)+d_{340}d_{10}\right\}\notag\\
&=\frac{3}{4p_1^2d_{10}^3}\left\{2p_1^3(q_2^2-q_1^2)+4p_1q_3(q_1+q_2)+d_{340}d_{10}\right\}.\label{Ef2c}
\end{align}
Finally, we have from (\ref{Qb}) and (\ref{Ef0}) that
\begin{equation}
\frac{1}{2}E(Q_n)=\frac{1}{2}+n^{-1}R+O(n^{-3/2}),\label{EQ}
\end{equation}
where with $\rho_1=n^{1/2}E(f_1)$ as in (\ref{f1}) and $\rho_2=E(f_2)$ as in (\ref{Ef2c}), we have
\begin{equation}
R=\rho_0+\rho_1+\rho_2,\label{rho}
\end{equation}
where $\rho_0$ is as in (\ref{rho0b}).
\newpage
This expression may be simplified into
$$R=-\frac{K}{24d_{10}^3p_1},$$
where
\begin{align*}
&K\\&=6 d_{10}^2 p_1 \{-2 (p_{11} + 6 p_{11} s_{0;11} + 2 s_1)
+ p_1 (2 + 12 s_1 + 5 s_1^2 - 17 s_2)\}\\
&+ 3\Big(4 p_{131}^2 - 8 p_1 p_{131} \{p_{21} + p_{11} (1 + s_{0;11} + s_{1;11}^{(12)})\}\\
&+ p_1^2 [8 p_{131} + p_{21}^2 + 24 p_{131} s_1 + 8 p_{131} s_1^2
+ 4 p_{11} p_{21} (5 + 5 s_{0;11} - 9 s_1 + 5 s_{1;11}^{(12)})\\
&- 4 p_{11}^2 \{15 s_1^2 - 18 s_1 (1 + s_{1;11}^{(12)}) + 2 (1 + s_{1;11}^{(12)})^2+ s_{0;11} (4 - 6 s_1 + 4 s_{1;11}^{(12)})\} - 8 p_{131} s_2]\\
&+ p_1^4 \{1 + 40 s_1^3 + 9 s_1^4 + 10 s_2 - 15 s_2^2 + 8 s_1 (3 + s_2)
+ s_1^2 (58 + 6 s_2)\}\\
&+ 2 p_1^3 [-p_{21} (7 + 12 s_1 + s_1^2 - 13 s_2)
+ 2 p_{11} \{-3 s_1^3 - s_1^2 (17 + 5 s_{1;11}^{(12)})\\
&- 3 s_1 (5 + 2 s_{1;11}^{(12)} - 5 s_2)
+ s_{0;11} (1 + 2 s_1 + 3 s_1^2 - 3 s_2)- (1 + s_{1;11}^{(12)}) (-1 +  7 s_2)\}]\Big)\\
&+ d_{10}\Big(-8 p_{131} - 6 p_{1761}
+ 3 p_1 (4 p_{11}^2 + 4 p_{21} - p_{241} + 4 p_{451} + 4 p_{131} s_{0;131}^{(11)} + 8 p_{11} s_1)\\
&+4 p_1^2 [-1 - 3 p_{21} - 6 s_1 + 27 p_{21} s_1 - 3 s_1^2 + 6 p_{12} s_1 s_{0;12}- 6 p_{12} s_{1;12}\\
& - 3 p_{21} s_{1;21}^{(43)} - s_2 - 6 p_{21} s_{0;21}
+ 3 p_{11} \{-1 + 2 s_{0;11} + 18 s_1^2 + 2 s_{1;11} + 2 s_{1;11}^{(12)}
+ (s_{1;11}^{(12)})^2\\
&- s_1 (19 + s_{1;11} + 16 s_{1;11}^{(12)})+ 3 s_2 + s_{2;11}^{(151)} - s_{2;121}^{(144)}\}]\\
&+ 3 p_1^3 \{3 - 4 (-9 + m) s_1^2 + 4 s_1^3 + 4 (-4 + m) s_2
- 4 s_1 (-7 + 4 s_2 + s_2^{(11)}) - 6 s_3 + 4 s_3^{(11)}\}\Big).
\end{align*}
Here, via (\ref{md10}) and (\ref{C10b}),
$$d_{10}=p_{11}-p_1(1+s_1).$$

In particular, for $m=2$, we get
$$R=\frac{1}{12\mu_1\mu_2}(1+9\mu_1+9\mu_2+15\mu_1\mu_2).$$

 \section{Bartlett correction}
In proposition 1 and corollary 2, asymptotically the expectation of $Q_n|\hat\lambda>0$ approaches 1, and we write $E(Q_\infty|\hat\lambda>0)=1$.

Now, consider a situation where no conditioning is involved. If, for the unconditional distribution, a result such as proposition 1 is present, the principle of Bartlett correction is to adjust the asymptotic statistic $Q_\infty$ according to
\begin{equation}
Q_n\approx\frac{E(Q_n)}{E(Q_\infty)}Q_\infty\approx(1+n^{-1}R)Q_\infty.\label{Bartlett}
\end{equation} 
It is clear from (\ref{Bartlett}) that the first moment of $\{E(Q_n)/E(Q_\infty)\}Q_\infty$ equals the first moment of $Q_n$. In fact, Lawley (1956) shows that under very general conditions, this is in fact true even for all higher moments.

Now, define the $\alpha$ percentile of $Q_n$, $q_n(\alpha)$, through
$\alpha=P(Q_n>q_n(\alpha))$. Then, from (\ref{Bartlett}),
$$\alpha\approx P\left\{(1+n^{-1}R)Q_\infty>q_n(\alpha)\right\}
=P\left\{Q_\infty>\frac{q_n(\alpha)}{1+n^{-1}R}\right\}.$$
Now, in the standard situation where $Q_\infty$ is $\chi^2_1$ distributed, we have\\ $\alpha=P\{Q_\infty>\chi^2_1(\alpha)\}$, i.e.
\begin{equation}
q_n(\alpha)\approx(1+n^{-1}R)\chi^2_1(\alpha).\label{qn}
\end{equation}
However, in our case we have that the distribution of $Q_\infty$ {\it conditional on} $\hat\lambda>0$ is $\chi^2_1$. With $\pi_n=P(\hat\lambda>0)$, we have
$$\alpha=P\{Q_n>q_n(\alpha)\}=\pi_nP\{Q_n>q_n(\alpha)|\hat\lambda>0\},$$
i.e. $\pi_n^{-1}\alpha=P\{Q_n>q_n(\alpha)|\hat\lambda>0\}$. Then, analogous to (\ref{qn}), we have
\begin{equation}
q_n(\alpha)\approx(1+n^{-1}R)\chi^2_1(\pi_n^{-1}\alpha).\label{qnb}
\end{equation}
Furthermore, we have that as $n\to\infty$, $\pi_n\to 1/2$. In a practical situation, we do not know $\pi_n$. Neither do we know $R$, but we can replace it by an estimate, $\hat R$, where we insert $\bar y_j$ for $\mu_j$, for all $j=1,2,...,m$. In particular, if $m=2$, corollary 2 gives us
\begin{equation}
\hat R=\frac{1}{6\bar y_1\bar y_2}(1+9\bar y_1+9\bar y_2+15\bar y_1\bar y_2).\label{Rhat}
\end{equation}
Hence, putting $\pi_n=1/2$, the practical analogue of (\ref{qnb}) is
\begin{equation}
\hat q_n(\alpha)=(1+n^{-1}\hat R)\chi^2_1(2\alpha).\label{qnhat}
\end{equation}

\section{Simulations}

To evaluate the Bartlett correction in a simulation study, we may compare estimates of
\begin{enumerate}
	\item $P\{Q_n>\chi^2_1(2\alpha)\}$
	\item $P\{Q_n>(1+n^{-1}\hat R)\chi^2_1(2\alpha)\}$
	\item $P\{Q_n>(1+n^{-1} R)\chi^2_1(2\alpha)\}$
	\item $P\{Q_n>(1+n^{-1} R)\chi^2_1(\hat\pi_n^{-1}\alpha)\}$
\end{enumerate}
These should go from worse to better. Note that the first one is the asymptotic approximation, the second one is the workable version of the Bartlett approximation while the other two are not useful in practice, but interesting by means of comparison. In the last one, $\hat\pi_n$ simply denotes the empirical proportion of positive $\hat\lambda$ in the simulation.

To check the validity and usefulness of the Bartlett correction, we have performed a small simulation study. Matlab R2109a was used. In each case, 1 000 000 replications were run. Optimization was performed using the procedure \texttt{fmincon}, and program codes are available upon request.

Below, in each case studied we depict expected simulated statistics together with the Bartlett approximation from proposition 1 (corollary 2 for $m=2$). Observe that asymptotically, $E(Q_n)$ tends to $1/2$, and in the figures we give $2E(Q_n)$ as the expectation. This is motivated by the fact that\\
$E(Q_n)=\pi_n*E(Q_n|\hat\lambda>0)$, i.e. $2E(Q_n)\approx E(Q_n|\hat\lambda>0)$ when $\pi_n\approx 1/2$. Because of this, we give the fractions of positive $\hat\lambda$ for different models in a separate graph at the end of the section. We also separately graph the simulated variances of $Q_n$ for the different models.

Casewise, we also graph the empirical rejection probabilities (\texttt{erp}) under the null  when using the asymptotic distribution, the Bartlett correction, the simulated Bartlett correction or (in some cases) the latter which, in addition, takes the simulated proportion of positive $\hat\lambda$ into account. In all cases, the nominal test sizes are 5\%.

We start by looking at dimension $m=2$. The first case studied is when $\mu_1=\mu_2=1$. The sample size goes from $n=20$ in steps of 20 up to $n=140$. We give the expected statistic in figure 1. We find that for small $n$, the Bartlett approximation overestimates the simulated expectation, but the accuracy improves as $n$ grows. Part of the explanation of the deviation is that the fraction of positive $\hat\lambda$, shown in figure 9, is lower than 0.5 for small samples, which means that the up weighting of the simulated $E(Q_n)$ by a factor 2 is too little in these cases.

As for the empirical test sizes (in per cent) in figure 2, we can see that the asymptotic \texttt{erp} values are too high for small $n$. This is no surprise, since the asymptotic expectation is too large, and also the simulated variances (see figure 10) are large. In view of figure 1, it is also to be expected that the Bartlett approximation over corrects, producing sizes a bit below the nominal 5\%. The simulated Bartlett correction does not suffer from this problem, and it also produces fairly accurate \texttt{erp} values.
 
It may perhaps be a bit surprising to see that, when taking positive $\hat\lambda$ proportions into account, we get slight over estimations. But this is probably due to the fact that the simulated variances are a bit large for small $n$, see figure 10. A large variance shifts the right tail of the distribution to the right, whereas a small fraction of positive $\hat\lambda$ shifts it to the left. In this case, when it comes to the simulated Bartlett correction with no correction for $\hat\lambda$, it seems that these two effects compensate each other. Imposing the $\hat\lambda$ correction, only the variance effect remains, and we get over estimation.

The cases $\mu_1=\mu_2=5$ and $\mu_1=\mu_2=20$ are plotted in figures  3-4 and 5-6, respectively. Here, the Bartlett corrections turn out to work much better. Observe that $\hat\lambda$ corrected values are not given in these figures, the reason being that $\hat\lambda$ is so close to 0.5 (see figure 9) so that these are virtually the same as the values obtained by the simulated Bartlett correction.

Regrettably, for dimensions $m>2$ the Bartlett correction does not seem so useful in practice. One example is the case with $m=3$ and\\ $\mu_1=\mu_2=\mu_3=1$, that is depicted in figures 7 and 8. As is seen from figure 7, the main problem is that, for the studied sample sizes $n$, the simulated expectation does not tend to one in a linear way as $n$ increases. (Supposedly, it does so for large enough $n$, but then the asymptotic distribution probably works so well that Bartlett correction is not relevant.)

Still, it seems from figure 8 that Bartlett correction works quite well for $n$ values of 100 and higher. But this is merely a stroke of luck, because as is seen from figures 9 and 10, the proportions of positive $\hat\lambda$ as well as the variances of $Q_n$ are low here, and this compensates for the fact that the Bartlett approximation over estimates the simulated expectation. This is also why the simulated Bartlett correction over estimates the nominal test size, even more so when imposing the $\hat\lambda$ correction.

Many more cases were simulated by the author, and further results may be obtained at request.

Finally, figure 9 gives the estimated proportions of positive $\hat\lambda$ for the models studied, and here we see that these are closer to 0.5 in cases with larger $\mu_j$ and further away in the dimension 3 case. In figure 10, we find that the variances of $Q_n$ are very similar in the dimension 2 cases, and a little smaller in the dimension 3 case. Observe that, as is easily seen, the asymptotic variance of $Q_n$ is $5/4=1.25$.

\section{Conclusion}

In the present paper, we have presented Bartlett corrections for the maximum log likelihood test of independence in the multivariate Poisson model proposed by e.g. Karlis (2003). This topic may be considered of interest, if not only by itself, so for giving an example of the Bartlett correction in the non standard type of situation where the parameter of the null hypothesis lies on the border of the parameter space.

In general, the resulting Bartlett correction factor is given by a very complicated formula, although explicit. But in the special case of dimension two (which is the case to be applied in Larsson, 2020), it simplifies into a nice and short expression. Our simulation studies also indicate that it is in dimension two that the correction is useful in practice. In higher dimensions, it seems that for sample sizes of practical interest, the expected statistic is far from a linear function of the sample size.

One possible improvement of the correction would be to incorporate the probability of a positive parameter estimate in some way. Another improvement could be to make the correction non linear. Computer intensive methods could also come into play. These would be interesting scopes for future research on the subject.

\section*{References}
Bartlett, M.S. (1937) Properties of sufficiency and statistical tests.
{\it Proceedings of the Royal Statistical Society, Series A}, 160, 268-282.
\vskip.2cm\noindent
Hayakawa, T. (1977) The likelihood ratio criterion and the asymptotic expansion of its distribution. {\it Annals of the Institute of Statistical Mathematics} 29, 359-378.
\vskip.2cm\noindent
Karlis, D. (2003) An EM algorithm for multivariate Poisson distribution and related models. {\it Journal of Applied Statistics} 30, 63-77.
\vskip.2cm\noindent
Larsson, R. (2020) Discrete factor analysis using a dependent Poisson model. {\it Computational Statistics} 35, 1133-1152.
\vskip.2cm\noindent
Lawley, D.N. (1956) A general method for approximating to the the distribution of the likelihood ratio criteria. {\it Biometrika} 43, 295-303.
\vskip.2cm\noindent
Self, S.G., Liang, K.Y. (1987) Asymptotic properties of maximum likelihood estimators and likelihood ratio tests under nonstandard conditions. {\it Journal of the American Statistical Association} 82, 605-610.

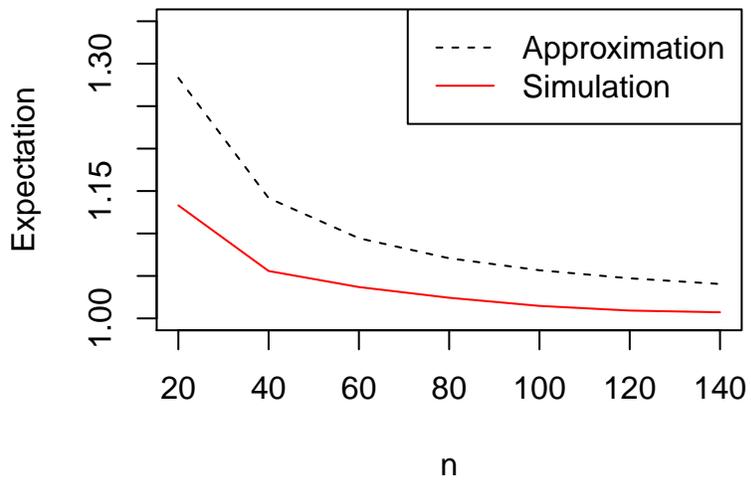
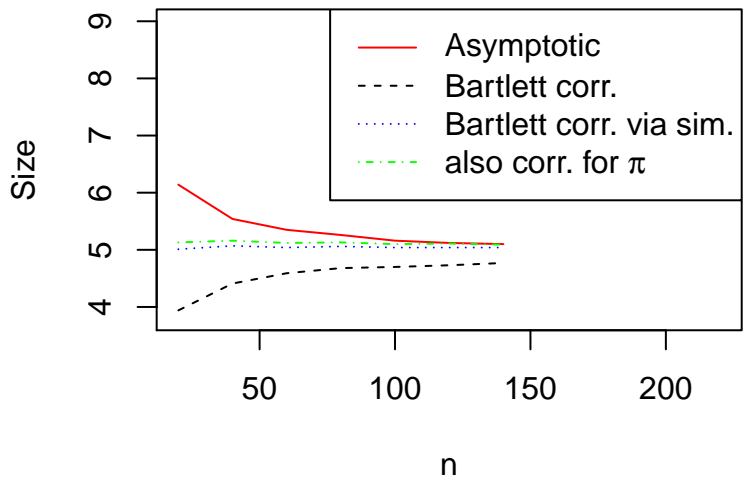
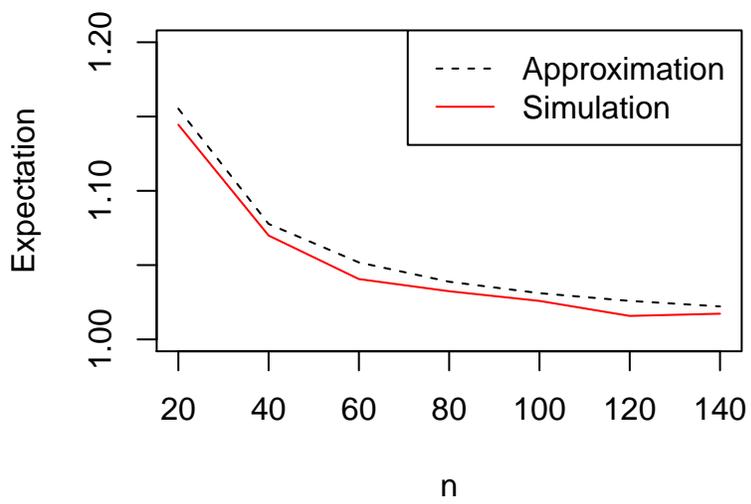
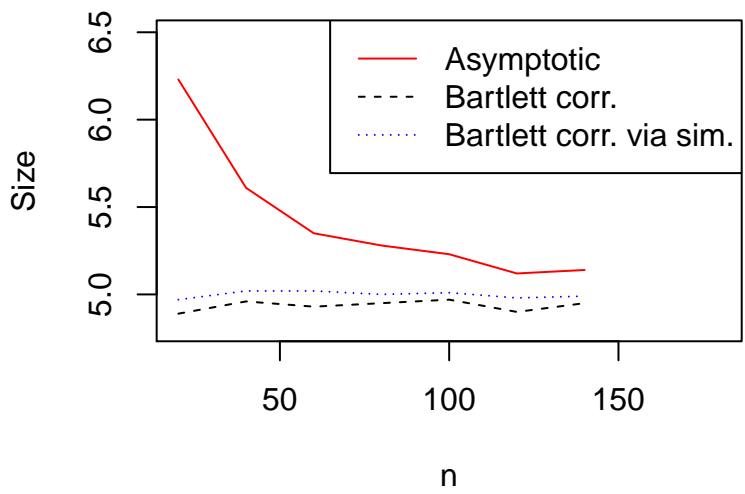
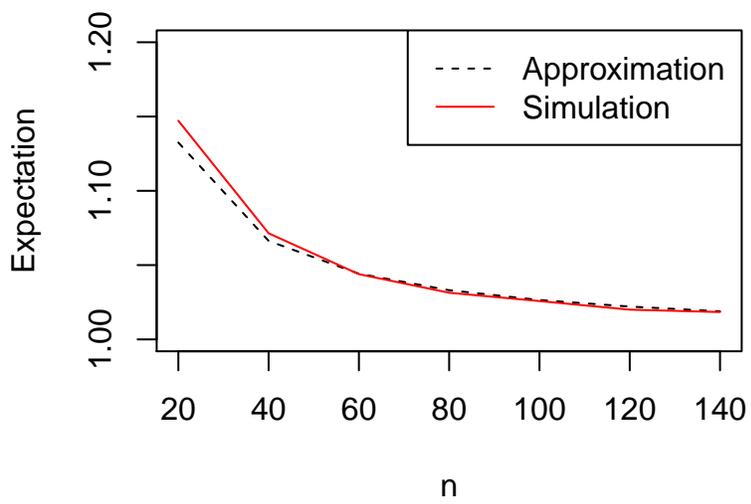
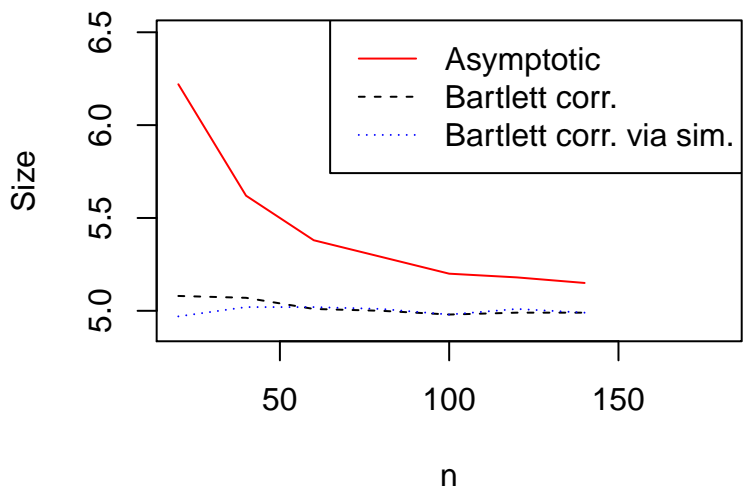
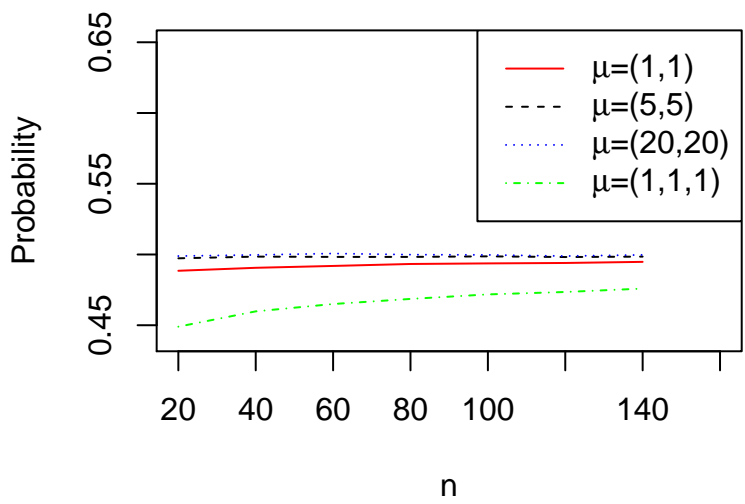
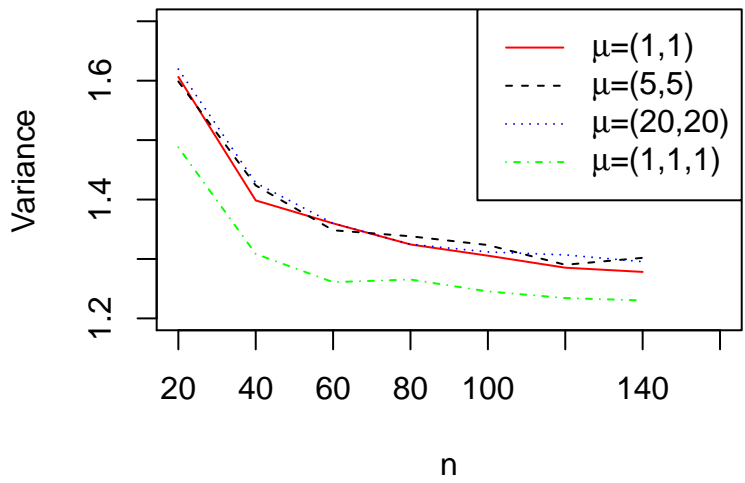
\end{document}